\documentclass[12pt]{amsart}
\usepackage{amsfonts}
\usepackage{epsfig,graphics}
\usepackage{amssymb}
\usepackage{amscd}

\newcommand{\BZ}{{\mathbb{Z}}}
\newcommand{\BN}{{\mathbb{N}}} 
\newcommand{\BR}{{\mathbb{R}}}
\newcommand{\BC}{{\mathbb{C}}}

\newcommand{\BQ}{{\mathbb{Q}}}

\newcommand{\BH}{{\mathbb{H}}}

\newcommand{\gD}{\Delta}
\newcommand{\gd}{\delta}
\newcommand{\gb}{\beta}

\newcommand{\gC}{\Gamma}
\newcommand{\gc}{\gamma}
\newcommand{\gs}{\sigma}
\newcommand{\gS}{\Sigma}
\newcommand{\gO}{\Omega}
\newcommand{\gep}{\epsilon}
\newcommand{\gl}{\lambda}
\newcommand{\ga}{\alpha}
\newcommand{\gt}{\tau}

\newcommand{\ti}[1]{\tilde{#1}}

\newcommand{\rk}{\text{rank}}
\newcommand{\rank}{\text{rank}}
\newcommand{\vol}{\text{vol}}

\newcommand{\SL}{\text{SL}}

\newcommand{\PSL}{\text{PSL}}

\newcommand{\SO}{\text{SO}}
\newcommand{\PSO}{\text{PSO}}

\newtheorem{prop}{Proposition}[section]
\newtheorem{thm}[prop]{Theorem}
\newtheorem{lem}[prop]{Lemma}
\newtheorem{cor}[prop]{Corollary}
\newtheorem{conj}[prop]{Conjecture}
\theoremstyle{definition}

\newtheorem{defn}[prop]{Definition}
\newtheorem{rem}[prop]{Remark}

\newtheorem{clm}[prop]{Claim}
\newtheorem{Ack}[prop]{Acknowledgments}

\begin{document}
\author{Tsachik Gelander}
\address{Tsachik Gelander, Institute of Mathematics, Hebrew university, Jerusalem, Israel}
\email{tsachik@math.huji.ac.il}
\thanks{Research partially supported by the Clore and the Marie Curie 
fellowships.}

\title{Homotopy type and volume of locally symmetric manifolds}
\maketitle

\begin{center}
\epsfig{file=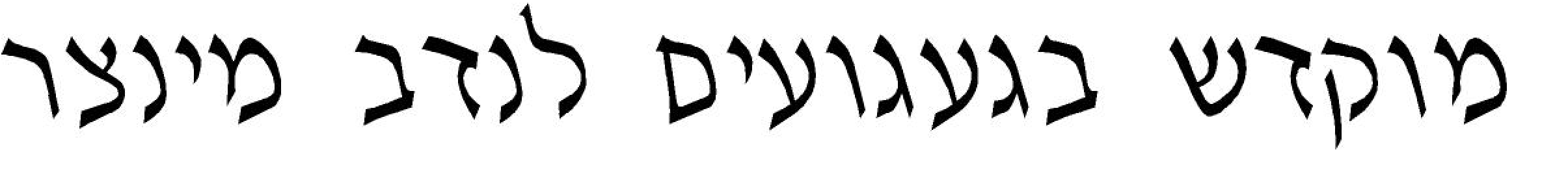,height=0.5cm}
\end{center}

\begin{abstract}
We consider locally symmetric manifolds with a fixed universal 
covering, and construct for each such manifold $M$ a simplicial complex 
$\mathcal{R}$ whose size is proportional to the volume of $M$. 
When $M$ is non-compact, $\mathcal{R}$ is homotopically equivalent to $M$, 
while when $M$ is compact, $\mathcal{R}$ is homotopically equivalent to 
$M\setminus N$, where $N$ is a finite union of submanifolds of fairly 
smaller dimension.
This reflects how the volume controls the topological structure of $M$, and 
yields concrete bounds for various finiteness statements which previously had 
no quantitative proofs. 
For example, it gives an explicit upper bound for the 
possible number of locally symmetric manifolds of volume bounded by $v>0$,
and it yields an estimate for the size of a minimal presentation for the 
fundamental group of a manifold in terms of its volume.
It also yields a number of new finiteness results.
\end{abstract}

\setcounter{tocdepth}{1}
\tableofcontents


\section{Introduction and statements of the main results}

In this article we study relations between the volume and the topological 
structure of locally symmetric manifolds. 
We are interested in asymptotic properties, 
when the volume tends to infinity. 

\medskip

We shall always fix a symmetric space, $S$, of non-compact type without 
Euclidean de-Rham factors, and consider the class of $S$-manifolds, 
by which we mean complete Riemannian manifolds locally isometric to $S$, 
or equivalently, manifolds of the form
$M=\gC\backslash S$ where $\gC$ is a discrete torsion-free group of isometries of $S$.
Sometimes we shall restrict our attention to arithmetic manifolds,
i.e. to Riemannian manifolds of the form 
$M=\gC\backslash S$, where $\gC\leq\textrm{Isom}(S)$ is a torsion-free arithmetic lattice
in the Lie group $\textrm{Isom}(S)$ of isometries of $S$.
A theorem of Borel and Harish-Chandra \cite{BoHC} says that if $M$ is arithmetic then 
$\vol (M)<\infty$. When $\rk (S)\geq 2$, Margulis' arithmeticity theorem \cite{Mar2}
gives the converse, i.e. $\vol (M)<\infty$ iff $M$ is arithmetic\footnote{The equivalence between finite
volume and arithmeticity is also known for the rank-1 cases $\text{Sp}(n,1)$ and $\text{F}_{4}^{-20}$ 
by \cite{GS}, \cite{cor}}. 

The ``complexity'' of the topology of locally symmetric manifolds is
controlled by the volume. This is illustrated by a theorem of 
Gromov (see \cite{BGS} theorem 2) which asserts  
that the Betti numbers are bounded by a constant times the volume, i.e. 
$$
\sum_{i=1}^{n}b_i(M)\leq c(S)\cdot\vol (M)
$$ 
for any $S$-manifold $M$.
Gromov's theorem applies also to non-locally symmetric manifolds, under appropriate conditions.

We conjecture (and prove in many cases) that, for locally symmetric manifolds,
the volume forces stronger topological restrictions. 

\begin{defn}\label{(d,v)-s.c.}
A {\bf $(d,v)$-simplicial complex} is a simplicial complex with 
at most $v$ vertices, all of them of valence $\leq d$.
\end{defn}

\begin{rem}
The number of $k$-simplexes in a $(d,v)$-simplicial complex $\mathcal{R}$ is 
$\leq\frac{v}{k+1}{n\choose k}$. Thus, the size of $\mathcal{R}$ (the number of its simplexes)
is at most $v\cdot\sum_{k=0}^d\frac{1}{k+1}{n\choose k}$, and this depends linearly on 
$v$.
\end{rem}          

\begin{conj}\label{conjA}
For any symmetric space of non-compact type $S$,
there are constants $\ga (S),d(S)$, such that any irreducible
$S$-manifold $M=\gC\backslash S$ (which is assumed also to be arithmetic in the case
$\dim (S)=3$) is homotopically equivalent to a 
$\big( d(S),\ga (S)\vol (M)\big)$-simplicial complex.
\end{conj}

\begin{rem}
The analogous statement is false for non-arithmetic manifolds in dimension $3$.
\end{rem}

In this paper we shall establish the following partial answers to Conjecture \ref{conjA}.
Since we failed in proving it in full generality, our results cannot be organized in a compact form, and
we have to split our statements and formulate strong results under appropriate conditions, and 
weaker results for more general cases.

\begin{thm}\label{thmA}
Let $S$ be a symmetric space of non-compact type. Then:
\begin{enumerate}
\item
Conjecture \ref{conjA} holds for non-compact arithmetic $S$-manifolds.
\item
If $S$ is neither isometric to $\SL_3(\BR )/\SO_3(\BR ),~\BH^2\times\BH^2$ nor to
$\BH^3$, then for some constants $\ga (S),d(S)$, the fundamental
group $\pi_1(M)$ of any 
$S$-manifold $M$ is isomorphic to the fundamental group of some
$\big( d(S),\ga (S)\vol (M)\big)$-simplicial complex.
\item
If $S$ is isometric to $\SL_3(\BR )/\SO_3(\BR ),~\BH^2\times\BH^2$ or to
$\BH^3$, then the fundamental group $\pi_1(M)$ of any  
$S$-manifold $M$ is a quotient of the fundamental group of some
$\big( d(S),\ga (S)\vol (M)\big)$-simplicial complex.
\end{enumerate}
\end{thm}

\begin{rem}
For compact arithmetic manifolds, Conjecture \ref{conjA} would follow, 
by a straightforward argument, if one could prove that the infimum of the
lengths of closed geodesics, taken over all compact arithmetic $S$-manifolds, 
is strictly positive. This conjectured phenomenon is strongly related to some 
properties of algebraic integers, such as the Lehmer conjecture,
which are still a mystery.
\end{rem}

Conjecture \ref{conjA} and Theorem \ref{thmA} yield
quantitative versions for some classical finiteness 
statements. We shall describe our main two applications 
in paragraphs \ref{BMP} and \ref{QVWT}.


\subsection{A linear bound on the size of a minimal presentation}\label{BMP}

\medskip

It is well known that the fundamental group of a locally symmetric manifold
with finite volume, or more generally a lattice in a connected semisimple Lie
group, is finitely presented. 
The compact case is quite standard, but the non-compact case was proved, step by step, by several authors
over several years.
Garland and Raghunathan \cite{Ga-Ra} proved it in the rank one case.
In the higher rank case, the finite generation was proved by Kazhdan \cite{Kazhdan} 
by defining and proving property-$T$ when $\rank (G)>2$,
and then by S.P. Wang \cite{S-P.W} when $\rank (G)=2$.
Using the finite generation, Margulis \cite{Mar2} proved the classical arithmeticity 
theorem for higher rank irreducible lattices\footnote{Margulis proved arithmeticity for higher rank 
non-uniform lattices \cite{Mar3}
few years before he proved the general case and without using his supper rigidity theorem}.
For arithmetic groups the finite presentability follows from the reduction 
theory of Borel and Harish-Chandra \cite{BoHC}.  
Later, using Morse theory, Gromov (\cite{BGS} theorem 2) gave a geometric proof of the finite 
presentability by proving that any locally 
symmetric manifold is diffeomorphic to the interior of a compact manifold with 
boundary. We remark that a slight modification of the argument of section \ref{6} below yields 
a completely elementary proof of this result of Gromov,
and hence of finite presentability (see also Remark \ref{finite-presentability}). 
However, in order to get a concrete estimate on the minimal possible size of such presentation 
in terms of the volume, we shall use arithmeticity (in the higher rank non-compact case).
We obtain the following quantitative version of finite presentability:

\begin{thm}\label{thmP}
Assume that $S$ is neither isometric to $\BH^3$,
$\SL_3(\BR )/\SO_3(\BR )$ nor to $\BH^2\times\BH^2$.
Then there is a constant $\eta =\eta (S)$ such that for any irreducible
$S$-manifold $M$, the fundamental group $\pi_1(M)$ admits a presentation
$$
 \pi_1(M)\cong\langle\Sigma :W\rangle
$$ 
with both $|\Sigma |,|W|\leq\eta\cdot\vol (M)$, in which all the relations 
$w\in W$ are of length $\leq 3$. 
\end{thm}

\begin{rem}
We believe (see also Conjecture \ref{conjA}) that the assumptions that $S$ is not
isometric to $\SL_3(\BR )/\SO_3(\BR )$ and to $\BH^2\times\BH^2$, are not really necessary.
However, our proof does not work in these cases.
It follows, however, from Theorem \ref{thmA}(2) that the analogous statements for non-compact 
$S$-manifolds hold also in these cases.
\end{rem}

\begin{rem}
The analogous statement for non-arithmetic hyperbolic $3$-manifolds is evidently false. 
However, if Conjecture \ref{conjA} is true, then the analogue of Theorem \ref{thmP} should
hold for arithmetic $3$-manifolds. By Theorem \ref{thmA}, it holds for non-compact arithmetic 
$3$-manifolds. Moreover, It was shown in \cite{cooper} that for any hyperbolic $3$-manifold $M$, 
the sum of the relations length, in any presentation of $\pi_1(M)$, is at 
least $\frac{\vol (M)}{\pi}$. This means that our upper bound is tight in 
this case.
\end{rem}

In the general case, Theorem \ref{thmA}(3) implies the following (weaker)
statement for which even the finiteness is in some sense surprising (in dimension $3$).
For a group $\gC$ let $d(\gC )$ denotes the minimal size of a generating set.

\begin{thm}\label{d(gC)}
For any $S$, there is a constant $\eta$ (depending on $S$), 
such that for any $S$-manifold $M$, $d\big(\pi_1(M)\big)\leq \eta\cdot\vol (M)$.
In other words, for any $v>0$, 
$$
 \sup\{d\big(\pi_1(M)\big) :\vol (M)\leq v\}\leq \eta v.
$$
Moreover, if $S$ is not isomorphic to $\SL_3(\BR )/\SO_3(\BR )$ then there is
a presentation
$$
\pi_1(M)\cong\langle\Sigma :W\rangle
$$ 
with $|\Sigma |,|W|\leq\eta\cdot\vol (M)$.
\end{thm} 

Note that Theorem \ref{d(gC)} does not give a bound on the length of the relations $w\in W$.


\subsection{A quantitative version of Wang's theorem}\label{QVWT}

\medskip

We apply our results in order to estimate the number of $S$-manifolds 
with bounded volume (or more generally the number of conjugacy classes of
lattices in $G=\text{Isom}(S)$). 
This can be considered as a continuous analogue to
asymptotic group theory which studies the subgroup growth
of discrete groups (where ``covolume of lattices'' extends the notion ``index of subgroups''). 
Unlike the situation in the discrete (finitely generated) case, even the 
finiteness statements are not clear, and in general not true. 
However, a classical theorem of H.C. Wang states that if $S$ is not isometric 
to one of the hyperbolic spaces $\BH^2,\BH^3$, then for any $v>0$ there 
are only finitely many irreducible $S$-manifolds with total volume $\leq v$, 
up to isometries (see \cite{Wa} 8.1, and paragraph \ref{SL_2} below). 
We remark that Wang's result and proof do not give explicit estimates.

Denote by $\rho_S(v)$ the number of non-isometric irreducible $S$-manifolds
with volume $\leq v$.

By Mostow's rigidity theorem, a locally symmetric manifold of dimension 
$\geq 3$ is determined by its fundamental group. Applying 
Theorem \ref{thmA}(2), we obtain: 

\begin{thm}
If $\dim (S)\geq 4$ and $S$ is neither isometric to $\SL_3(\BR )/\SO_3(\BR )$ 
nor to $\BH^2\times\BH^2$ then there is a constant $c$, depending on $S$, with
respect to which
$$
\log\rho_S(v)\leq c\cdot v\log v
$$ 
for any $v>0$.
\end{thm}

This upper bound was first proved for hyperbolic manifolds by M. Burger, A. 
Lubotzky, S. Mozes and the author in \cite{BGLM}, where also a lower bound of
the same type was established, proving that this estimate is the true 
asymptotic behavior in the hyperbolic case. 

\begin{thm}[BGLM]
For $n\geq 4$, there are constants $c_n>b_n>0$ and $v_n>0$ such that
$$
 b_nv\log v\leq \log\rho_{\BH^n}(v)\leq c_nv\log v,
$$
whenever $v>v_n$.
\end{thm}

However, we suspect that in the higher rank case, where all manifolds are
arithmetic and conjectured to possess the congruence subgroup property, 
this upper bound is far from the true asymptotic behavior.  
It seems, in view of \cite{Lub}, that the problem of determining the real
asymptotic behavior is closely related to the congruence subgroup problem.
It is also related to the analysis of the Galois cohomology of compact 
extensions of $G=\text{Isom}(S)$.

For hyperbolic manifolds, the weaker upper bound  
$\rho_{\BH^n}(v)\leq v\exp (\exp (\exp (v+n)))$ was proved previously by Gromov \cite{Gr1}.
In this paper we establish a first concrete estimate for $\rho_S(v)$ for general $S$.

In dimension $2$ and $3$ the analogue of Wang's finiteness theorem is false.
However, as was shown by Borel \cite{Bor}, it remains true when 
considering only arithmetic manifolds.
The following estimate for the number of non-compact arithmetic $3$-manifolds
follows from Theorem \ref{thmA}(1):

\begin{prop}\label{116}
For some constant $c>0$, there are at most $v^{cv}$ non-isometric arithmetic 
non-compact hyperbolic 3-manifolds with volume $\leq v$.
\end{prop}

For compact arithmetic $3$-manifolds we conjecture that the analogous 
statement holds, but prove only the following weaker statement:

\begin{prop}\label{117}
Let $M$ be a compact (arithmetic) hyperbolic 3-manifold.
Then for some constant $c(M)$, the number of non-isometric hyperbolic
manifolds commensurable to $M$, with volume $\leq v$, is at most $v^{c(M)v}$.
\end{prop}

We remark that the analogue of Propositions \ref{116} and \ref{117} hold also in 
the cases $S=\SL_3(\BR )/\SO_3(\BR )$ and $S=\BH^2\times\BH^2$.

\medskip

In section \ref{complements}, we shall generalize some of our results concerning $S$-manifolds
to the larger family of $S$-orbifolds. We shall also indicate how to construct triangulations
for rank-1 manifolds which are not necessarily arithmetic.

\medskip

Let us now give a short and not precise explanation of the basic lines of the proofs of our main 
results \ref{thmA}.
The idea is to construct, inside each $S$-manifold $M$, a submanifold with boundary $M'$, which is
similar enough to $M$ and for which we can construct a triangulation of size $\leq c\cdot\vol (M)$.
In order to construct such a triangulation we shall require a lower bound on the injectivity radius
of $M'$ (independent of $M$) and some bounds on the geometry of the boundary $\partial M'$.
We shall construct $M'$ inside the $\gep_s$-thick part, in a way that its pre-image in the 
universal cover $S$ will be the complement of some ``locally finite'' union of convex sets, each
has a smooth boundary whose curvature is bounded uniformly from below, and the angles at the 
corners of $\partial M'$ (where the boundaries of two or more such convex sets meet) are bounded 
uniformly from below.

The non-compact arithmetic case is easier to deal with. The reason is that any non-uniform 
arithmetic lattice $\gC\leq G=\text{Isom}(S)^0$ comes from a rational structure on $G$ (rather then
on a compact extension of $G$ as the case might be when $\gC$ is uniform). 
This implies that there is 
some constant $\gep_s'$ such that any non-uniform arithmetic $S$-manifold contains no closed 
geodesics of length $\leq\gep_s'$ (see section \ref{AML}). In other words, any non-trivial closed 
loop which is short enough corresponds to a unipotent element in the fundamental group.
We define (the pre-image in $S$ of) $M'$ to be the complement of the union of appropriate
sub-level sets for the displacement functions $\{ d_\gc\}$ where $\gc$ runs over all unipotents
in $\pi_1(M)$. The injectivity radius in $M'$ is then uniformly bounded from below,
and using the fact that unipotents acts nicely on $S$ and on $S(\infty )$ we can 
both estimate the 
geometry of $\partial M'$, and also construct a deformation retract from $M$ to $M'$.

The situation is more subtle in the compact case. We do not have sufficient information
on the thick-thin decomposition. 
In this case we construct $M'$ and prove that it is diffeomorphic to $M\setminus N$ where $N$ is a 
finite union of submanifolds of codimension $\geq 3$. $N$ will contain the union of all 
closed geodesics of length smaller than some fixed constant. 
The idea is that at any point outside $N$
there is a preferred direction, such that when we move along it, the injectivity radius is
increasing most rapidly. This will help us to define a deformation retract from $M\setminus N$ to 
$M'$. As $N$ has codimension $\geq 3$, $\pi_1(M)\cong\pi_1(M\setminus N)\cong\pi_1(M ')$.
A main difficulty (which arises also in the compact rank one case) is how to control the geometry 
of the boundary of $M'$. We shall handle this difficulty in section \ref{bigangles}, where the 
main idea is Lemma \ref{bigA} which says that if two isometries commute then the exterior angle 
between their sub-level sets is $\geq\frac{\pi}{2}$.

\medskip

The present paper generalizes the work \cite{BGLM} which treated the special case of real 
hyperbolic spaces. Although some of the ideas from \cite{BGLM} appear again in the present paper, 
the situation in the general case is significantly more complicated than the hyperbolic case. 
The proof in \cite{BGLM} uses the explicit description of a hyperbolic compact 
thin component of the thick-thin decomposition as a cone over a coaxial Euclidean ellipsoids, 
as well as some computations in constant curvature. Hence the
argument in \cite{BGLM} does not apply to more general rank-1 symmetric spaces.
More crucially, in contrast to the situation in the rank-1 case, in the higher rank case
there is currently no good understanding of the structure of the thin components in the 
thick-thin decomposition, and hence new ingredients are required also in the skeleton of the proof.

Most of this paper is devoted to the treatment of the higher rank case. 
However, whenever it is not required, we shall not make
any assumption on the rank. We remark also that many of the arguments in this paper stay valid for
general manifolds of non positive curvature. In particular, some parts of this work could be 
generalized to non symmetric Hadamard spaces.


\section{Notations, definitions and background}
In this section we shall fix our notations and
summarize some basic facts about semisimple Lie groups, symmetric spaces of 
non-compact type, and manifolds of non-positive curvature. For a comprehensive
treatment of these subjects we refer the reader to \cite{Rag}, \cite{BGS} and \cite{BH}.  

Let $S$ be a symmetric space of non-compact type. We shall always assume, 
that $S$ has no Euclidean 
de Rham factors. Let $G=\textrm{Isom}(S)$ be the Lie group of isometries of $S$. $G$ is 
center-free, semi-simple without compact factors, and with finitely many connected components. 
We denote by $G^0$ the identity component of $G$ with respect to the real 
topology. There is an algebraic group $\Bbb{G}$ defined over $\BQ$, such that
$G^0$ coincides with the connected component of the group of real points $\Bbb{G}(\BR )^0$.
In particular $G^0$ admits, apart from the real topology, a Zariski topology
which is defined as the trace in $G^0$ of the Zariski topology in $\Bbb{G}$.
$G^0$ acts transitively on $S$, and  we can identify $S$ with $G^0/K$ where $K\leq G^0$ is
a maximal compact subgroup. We remark that there is a bijection between
symmetric spaces of non-compact type and connected center free semisimple Lie 
groups without compact factors.

$S$ is a Riemannian manifold with non-positive
curvature such that for each point $p\in S$ there is an isometry $\sigma_p$
of $S$ which stabilizes $p$ and whose differential $d_p(\sigma_p)$ at $p$ is $-1$. 
The composition of two such isometries $\sigma_p\cdot \sigma_q$ is called a 
{\bf transvection}, and is belong to $G^0$.
The non-positivity of the sectional curvature means that the distance 
function $d:S\times S\to \BR ^+$ is convex and, in fact, its restriction to 
a geodesic line $c(t)=\big( c_1(t),c_2(t)\big)$ in $S\times S$ is strictly 
convex, unless $c_1$ and $c_2$ are contained in a flat plane in $S$.

For a real valued function $f:X\to\BR$ on a set $X$ we denote by
$\{ f<t\}$ the $t$-sub-level set
$$
 \{ f<t\}=\{ x\in X:f(x)<t\}.
$$

For $\gc \in G$ we denote by $d_{\gc}$ the {\bf displacement function}
$$
  d_{\gc}(x)=d(\gc \cdot x,x).
$$
This function is convex and smooth outside Fix$(\gc )$. 
In particular the sub-level sets $\{ d_{\gc} <t\}$ are convex with smooth boundary.
We denote by {\bf min}$(\gc )$ the set 
$$
 \min (\gc )=\{x\in S:d_{\gc}(x)=\inf d_{\gc}\}.
$$

For a set $A\subset S$ we denote by 
$$
 D_A(x)=d(A,x)
$$ 
the {\bf distance function}.
If $A$ is a convex set then the function $D_A$ is convex and smooth at any
point outside the boundary of $A$. We denote its $t$-sub-level set by
$$
 (A)_t=\{ D_A<t\} 
$$
and call it the $t$-{\bf neighborhood} of $A$. 
Note that $(A)_{t+s}=\big( (A)_t\big)_s$.
Similarly, we let $)A(_t$ denote the $t$-{\bf shrinking} of $A$
$$
 )A(_t=S\setminus (S\setminus A)_t.
$$
If $A$ is convex then $\big) (A)_t\big(_t= A$ but in general these are 
different sets. 

A subset $\mathcal{N}$ of a metric space $X$ is called an {\bf $\gep$-net}, if
$(\mathcal{N})_\gep =X$. A subset $\mathcal{N}\subset X$ is called 
{\bf $\gep$-discrete} is $d(x,y)>\gep$ for any $x,y\in\mathcal{N},~x\neq y$.
If $\mathcal{N}$ is a maximal $\gep$-discrete subset, then it is an $\gep$-net.

If $A$ is closed and convex then for any $x\in S$ there is a unique closest 
point $p_A(x)$ in $A$. The map $p_A:S\to A$ is called the {\bf projection} 
on $A$, and it is distance decreasing, i.e.
$d\big(p_A(x),p_A(y)\big)\leq d(x,y)$ for any $x,y\in S$.

A {\bf flat subspace} is a totally geodesic sub-manifold which is 
isometric to a Euclidean space. A {\bf flat} is a maximal flat subspace. 
Any flat subspace is contained in a flat. 
All flats in $S$ have the same dimension
$\textrm{rank}(S)=\textrm{rank}(G)$\footnote{By $\textrm{rank}(G)$ we shall 
always mean the real rank of $G$}. A geodesic $c\subset S$ is called 
{\bf regular} if it is contained in a unique flat.

An element $\gc \in G$ is unipotent if $\text{Ad}(\gc )-1$ is a nilpotent 
endomorphism. A subgroup of $G$ is called a {\bf unipotent subgroup} if all its 
elements are unipotent. If $H\leq G^0$ is a unipotent subgroup then its 
Zariski closure $\overline{H}^z\leq\Bbb{G}$ is also unipotent, and 
$\overline{H}^z_\BR$ is connected. Moreover, a unipotent subgroup $H\leq G^0$ 
is connected iff $H=\overline{H}^z_\BR$. 

$S(\infty )$ is the set of equivalence classes of geodesic rays, where two
rays are considered equivalent if their traces are of bounded Hausdorff 
distance from each other. 
There is a structure of a spherical building on $S(\infty )$ - the Tits 
building of $S$. The chambers of this building are sometimes called Weyl 
chambers. A maximal unipotent subgroup of $G^0$ is the unipotent radical of 
some minimal parabolic subgroup. The minimal parabolic subgroups of $G^0$ are 
the stabilizers (and the fixator) of the chambers of this building, and there
are canonical bijections between the set of minimal parabolic subgroups, the 
set of maximal unipotent subgroups and the set of chambers. All maximal 
unipotent subgroups are conjugate in $G^0$, or in other words, $G^0$ acts 
transitively on the set of chambers $W\subset S(\infty )$. Moreover, each 
parabolic subgroup acts transitively on $S$, and hence $G$
acts transitively on the set of couples $(W,x)$ of a chamber 
$W\subset S(\infty )$ and a point $x\in S$.
There is a canonical metric on $S(\infty )$, called the Tits metric, with 
respect to which the apartments of the Tits building are isometric to spheres
and the chambers are isometric to spherical simplices. The induced action of 
$G$ on $S(\infty )$ preserves the Tits metric.
A point $p\in S(\infty )$ is called {\bf regular} if it is an interior point
of a chamber. For any $p\in S(\infty )$ and $x\in S$ there is a unique geodesic
$c$ with $c(0)=x$ and $c(\infty )=p$, and $p$ is regular iff $c$ is regular.

When $\gC$ is a group of isometries of $S$, we say that a subset $A\subset S$ 
is $\gC$-{\bf precisely invariant} if $\gc \cdot A=A$
whenever $\gc \cdot A \cap A\neq \emptyset$.   
If $\gC$ acts freely and discretely (i.e. $\gC \subset G$ is discrete
and torsion free) then we denote by $\gC\backslash A$ the image of $A$ in $\gC\backslash S$.
If $A$ is a connected simply connected $\gC$-precisely invariant set
then 
$$
 \pi_1(\gC\backslash A )\cong \{ \gc \in \gC :
 \gc \cdot A=A\}. 
$$

For a subset $M'$ of $M=\gC\backslash S$ we let $\ti M'$ denote its pre-image in $S$, and
we shall usually denote by $\ti{M'}^0$ an arbitrary connected component of 
$\ti M'$.  

For a subgroup $\gC\leq G$, a constant $\gep >0$ and a point $x\in S$,  
we let $\gC_\gep (x)$ denote the group generated by the ``small elements'' 
at $x$
$$
 \gC_\gep (x)=\langle \gc\in\gC :d_\gc (x)\leq\gep\rangle.
$$
Recall the classical Margulis lemma:

\begin{thm}[{\bf The Margulis lemma}]
There are constants $\gep_s>0$ and $n_s\in\BN$, depending only on $S$,
such that for any discrete subgroup $\gC\leq G$, and any point $x\in S$,
the group $\gC_{\gep_s}(x)$ contains a subgroup 
of index $\leq n_s$ which is contained in a connected nilpotent subgroup of 
$G$.
\end{thm}

The {\bf $\gep$-thick thin} decomposition of an $S$-manifold $M=\gC\backslash S$ reads 
$$
 M=M_{\geq\gep}\cup M_{\leq\gep}.
$$
The {\bf $\gep$-thick part} $M_{\geq\gep}$ is defined as the set of all 
points $x\in M$ at which the injectivity radius is $\geq\frac{\gep}{2}$, 
and the {\bf $\gep$-thin part} $M_{\leq\gep}$
is the complement of the interior of the $\gep$-thick part.
Note that $M_{\leq\gep}$ is the set of points in $M$ through which there is a
non-contractible closed loop of length $\leq\gep$.
The pre-image of $M_{\leq\gep}$ in $S$ has a nice description as a 
locally finite union of convex sets
$$
 \ti M_{\leq\gep}=\cup_{\gc\in\gC\setminus\{ 1\}}\{ d_\gc\leq\gep\}.
$$
The Margulis lemma yields an important piece of information on the structure 
of the $\gep$-thick-thin decomposition when $\gep\leq\gep_s$.

If $\gC \subset G$ is a uniform lattice then any element 
$\gc \in \gC$ is semisimple, i.e. $\text{Ad}(\gc )$ is diagonalizable over 
$\BC$, or equivalently $\min (\gc )\neq \emptyset$. If, in addition, $\gC$ is 
torsion free then all its elements are hyperbolic.
A semisimple element is hyperbolic iff it has a complex eigenvalue outside the
unit disk.
A hyperbolic isometry $\gc$ has an axis, i.e. a geodesic line on which
$\gc$ acts by translation. Any two axes  are parallel, and 
$x\in \min (\gc )$ iff the geodesic line $\overline{x,\gc \cdot x}$
is an axis of $\gc$. 
We define a {\bf regular isometry} to be a hyperbolic isometry whose axes are regular 
geodesics (if one axis is a regular geodesic, then so are all of them, since axes are parallel). 
If $g\in G$ acts
as a regular isometry then $g$ is a regular element in $G$ in the usual sense,
but not necessarily vice versa.
If $\gc\in G$ acts as a regular isometry, and $F$ is the unique flat containing
an axis of $\gc$ then $\gc$ preserves $F$ and acts by translation on it, and
in particular $\min (\gc )=F$.
(This follows for example from the facts that the Weyl group 
$N_{G^0}(A)/A$, for $A=C_{G^0}(\gc )$, acts simply transitively on the set of Weyl chambers 
in $F(\infty )$, and that a
regular point $p\in F(\infty )$ determines its ambient Weyl chamber.)

If $A\subset S$ is a closed convex set, and $\ga\in G$ preserves $A$, i.e.
$\ga\cdot A=A$, then the projection $P_A:S\to A$ commutes with $\ga$, and 
since $P_A$ does not increase distances, we have
$d_\ga \big( P_A(x)\big)\leq d_\ga(x)$. In particular, $\ga$ is semisimple iff
$\min (\ga )\cap A\neq\emptyset$. Since for any isometry $i$ of a flat 
subspace $F$, $\min (i)\neq\emptyset$, it follows that if $\ga\in G$ preserves
a flat subspace then $\ga$ is semisimple. If $g\in G^0$ is semisimple, then
its centralizer group $C_{G^0}(g)$ is a closed reductive group.
If $g\in G$ is a transvection with axis $c$, then
$$
 C_{G^0}(g)=\{ h\in G^0:h\cdot c \text{ is parallel to } c\}.
$$
   
If $\gC \subset G$ is a non-uniform lattice, then it is almost generated by
unipotents. $\gc$ is unipotent iff $\text{inf}(d_\gc )=0$ while $d_\gc$ never 
vanishes on $S$. In particular a unipotent element stabilizes a point at infinity but not in $S$.

If $f:G^0\to H$ is a surjective Lie homomorphism, and $g\in G$ is semisimple
(resp. unipotent) then $f(g)$ is semisimple (resp. unipotent).


\section{Some deformation retracts}\label{3}

In several places in this paper we will need to produce a deformation 
retract from a submanifold (maybe with boundary) to a subset of it which
is defined by the condition that all the functions from some given family
are larger than a given constant. A typical example of such a situation is 
given by the manifold itself and the $\gep$-thick part.
The family of functions in this example is indexed by the fundamental group.
To each homotopy class of closed loops corresponds the function for which 
the value at any point  $x$ is the minimal length of a loop from the class, 
which passes through $x$. The $\gep$-thick part is exactly the subset where 
all these functions are $\geq\gep$. The Margulis lemma yields information 
about the set of functions which are $<\gep$ at $x$, which can sometimes be exploited in 
order to define a deformation retract to the $\gep$-thick part.
Another example is given by a subset $M'$
of the manifold $M$ and its $\gep$-shrinking $)M'(_\gep$. If the complement of $M'$ is 
given by the union of subsets satisfying certain properties,
it is sometimes convenient to use the distance functions from these subsets
in order to define a deformation retract to $)M'(_\gep$.

The idea is to construct a continuous vector
field which makes an acute angle at any point $x$ with the gradients of all the
functions from the given family which are small at $x$, and then to let points
flow along its integral curves.
This idea has been used in an elegant way in other places (c.f. \cite{BGS}). The situation concerned here, however,
is more subtle since we consider manifolds with boundary, and have to make sure that on the 
boundary, the vector field is pointing towards the interior.

\begin{defn}
For a family of real valued continuous functions 
$\mathcal{F} =\{\phi_i\}_{i\in I}$ on
a manifold $Y$ (with or without boundary) and $\gep >0$, we define the {\bf  
$(\mathcal{F},\gep )$-thick part} (or simply the {\bf $\mathcal{F}$-thick part) 
$\mathcal{F}_{\geq\gep}$} as follows
$$
 \mathcal{F}_{\geq\gep}=\cap_{\phi\in\mathcal{F}}\overline{\{\phi >\gep\}}.
$$
When $\mathcal{F}$ is given, we denote by 
$\Psi_{x,\gt }~(x\in Y,\gt\in\BR )$ the subset
$$
 \Psi_{x,\gt } =\{\phi\in\mathcal{F} : \phi (x)\leq\gt\}.
$$
Abusing notations, we shall sometimes write $\Psi_{x,\gt }$ for the 
corresponding set of indices
$\Psi_{x,\gt } =\{i\in I: \phi_i(x)\leq\gt\}$.
We say that the family $\mathcal{F}$ is {\bf locally finite} if
$\Psi_{x,\gt }$ is finite for any $x\in Y$ and $\gt\in\BR$. 
We say that $\gt\in\BR$ is a {\bf critical value} of a continuous function
$\phi$ if $\overline{\{\phi >\gt\}}\neq \{\phi\geq\gt\}$, i.e. if $\gt$ is 
a value of a local maximum of $\phi$.
\end{defn}

If $Z\subset Y$ is a submanifold with smooth boundary, and $\phi$ a real 
valued function on $Z$. The gradient $\nabla\phi$ is well defined on $Z$ 
wherever it is unique (on 
the boundary it is the tangent vector with length and direction equal to the 
value of the maximal directional 
derivative and the direction at which it occurs). The function $\phi$ is said 
to be smooth (or $C^1$) if its gradient $\nabla\phi$ is a continuous function on the 
whole of $Z$. When $\phi$ is smooth, its directional derivative with respect to
$v\in T_x(Y)$ is given by the usual formula $v\cdot\nabla\phi (x)$.
 
In later sections, we will consider families $\mathcal{F}$ which consist of functions of the 
following 2 types:
\begin{itemize}
\item 
The displacement function $d_\gc (x)=d(x,\gc\cdot x)$ associated with
a non-trivial isometry $\gc$ of $S$.
\item 
The distance function $D_A(x)=d(x,A)$ from a given closed convex set with
smooth boundary.
\end{itemize}
In both cases the function is convex, non-negative, and without positive 
critical values. The function $d_\gc$ is smooth on $\{ d_\gc >0\}$, and the 
function $D_A$ is smooth on $\overline{\{ D_A>0\}}$.

\medskip


\begin{lem}[\bf A deformation retract which increases functions]
\label{deformation-retract}
Consider a submanifold without boundary $Y\subset S$.
\begin{itemize}
\item 
Let $\mathcal{F} =\{\phi_i\}_{i\in I}$ be a locally finite
family of non-negative continuous functions on $Y$.
\item
Assume that for each $\phi\in\mathcal{F}$ with $\{\phi =0\}\neq\emptyset$, the set 
$\overline{\{\phi >0\}}$ is a submanifold with smooth boundary.
\item
Assume that each $\phi\in\mathcal{F}$ is $C^1$ on $\overline{\{\phi >0\}}$.
\item 
Let $L=\mathcal{F}_{\geq 0}=\cap_{\phi\in\mathcal{F}}\overline{\{\phi >0\}}$.
\item
Let $\gb :[0,3\gep ]\to\BR^{>0}$ be a given continuous positive function 
(in case all $\phi\in\mathcal{F}$ are strictly positive we allow $\gb$ to
be defined only on $(0,3\gep ]$).
\item
Assume that $\gep$ is not a critical value for any $\phi\in\mathcal{F}$.
\item
Assume that for any $x\in Y$    
there is a unit tangent vector $\hat n(x)\in T_x(S)$ such that  
$$
 \hat n(x)\cdot \nabla \phi (x)\geq\gb\big( \phi (x)\big)
$$
for any $\phi\in\Psi_{x,3\gep}$. 
\end{itemize}
Then there is a deformation retract from $L$ to the $\mathcal{F}$-thick 
part $\mathcal{F}_{\geq\gep}$. 
(In particular it follows that $\mathcal{F}_{\geq\gep}\neq\emptyset$.)

If $\gC\leq\text{Isom}(S)$ is a discrete subgroup, and if $L$ and the family 
$\mathcal{F}$ are $\gC$-invariant, in the sense that $\gc\cdot L=L$, and
the function $\gc\cdot\phi (x):=\phi(\gc^{-1}\cdot x)$ belongs to $\mathcal{F}$ for any
$\phi\in\mathcal{F},~\gc\in\gC$,
then there exists such a deformation retract
which is $\gC$-invariant (in the obvious sense)
and hence induces a deformation retract between 
the corresponding subsets $\gC\backslash L$ and 
$\gC\backslash\mathcal{F}_{\geq\gep}$ of $M=\gC\backslash S$.
\end{lem}

\begin{proof}
We will define an appropriate continuous vector field on $L$. 
The desired deformation retract will be the flow along this vector field.

Let $\gd (x)$ denote 
$$
 \gd (x)=\min_{\phi\in \mathcal{F}}\phi (x).
$$

For any non-empty subset $\Psi\subset\mathcal{F}$ for which 
$\cap_{\phi\in\Psi}\{\phi\leq 3\gep\}\neq\emptyset$ and any point $x$ of this intersection, let 
$\hat f(x,\Psi )\in T_x(S)$ be a unit tangent vector which maximizes the 
expression
$$
 \min\{\hat f\cdot \nabla\phi (x):\phi\in\Psi\},
$$
and let 
$$
 \gb_\Psi (x)=\min_{\phi\in\Psi}\gb\big(\phi (x)\big).
$$
Then for any $\phi\in\Psi$
\begin{eqnarray*}
 \hat f(x,\Psi )\cdot \nabla \phi (x)\geq
 \min_{\phi\in\Psi}\hat f(x,\Psi )\cdot \nabla \phi (x)\geq
 \min_{\phi\in\Psi}\hat n(x)\cdot \nabla \phi (x)\geq
 \min_{\phi\in\Psi}\gb\big(\phi (x)\big)=
 \gb_\Psi (x).
\end{eqnarray*} 
Moreover, it follows from the strict convexity of the Euclidean unit disk that
$\hat f(x,\Psi )$ is uniquely determined, and consequently, that for a fixed
$\Psi$, the vector field $\hat f(x,\Psi )$ is continuous on
$\cap_{\phi\in\Psi}\{\phi\leq 3\gep\}$.

The desired vector field is defined as follows:
\begin{eqnarray*}
 &&\!\!\!\!\!\!\! \overrightarrow V(x)=
 \sqrt{2\big(\gep -\gd(x)\big)\lor 0}\cdot\\
&\phantom{\leq}&
 \cdot\sum_\Psi
 \frac{\big( 3\gep -\max_{\phi\in\Psi}\phi (x) \big)\lor 0}{\gep}\cdot
 \frac{\big(\min_{\phi\notin\Psi}\phi (x)-\gep \big)\lor 0}{\gep}\cdot
 \frac{1}{\gb_\Psi(x)}\hat f(x,\Psi ),
\end{eqnarray*}
where the sum is taken over all non-empty finite subsets 
$\Psi \subset\mathcal{F}$.

Notice that $\gb_\Psi(x)$ is strictly positive, and 
all the terms in each summand are continuous, and $\overrightarrow V\equiv 0$
on the $\mathcal{F}$-thick part
$$
 \mathcal{F}_{\geq\gep}=
 \cap_{\phi\in\mathcal{F}}\overline{\{\phi >\gep\}}=
 \cap_{\phi\in\mathcal{F}}\{\phi\geq\gep\}=
 \{\gd\geq\gep\}.
$$
The term $\sqrt{2\big(\gep -\gd(x)\big)\lor 0}$ takes care of the continuity
on the boundary $\{\gd =\gep\} =\partial\{\gd \leq\gep\}$ of the 
$\mathcal{F}$-thin part. 
The terms $\frac{\big( 3\gep -\max_{\phi\in\Psi}\phi (x) \big)\lor 0}{\gep}$ 
guarantee that all the non-zero summands correspond to sets which are 
contained in the finite set $\Psi_{x,3\gep}$. In particular the summation is 
finite for any $x\in L$, and $\hat f(x,\Psi )$ is well defined for any 
non-zero summand.
The terms $\frac{\big(\min_{\phi\notin\Psi}\phi (x)-\gep \big)\lor 0}{\gep}$ 
guarantee that all the non-zero summands correspond to $\Psi$'s which contains
$\Psi_{x,\gep}$.   

If $\gd (x)<\gep$ then $\Psi_{x,\gep}\neq\emptyset$, and, when looking only on the summand 
corresponding to $\Psi_{x,2\gep}$, we see that
$$
 \nabla \phi (x)\cdot\overrightarrow V(x)\geq\sqrt{2\big(\gep -\gd(x)\big)}
$$
for any $\phi\in\Psi_{x,\gep}\subset\Psi_{x,2\gep}$, and in particular for any 
$\phi$ with $\phi (x)=\gd (x)$.

It follows that if $x(t)$ is an integral curve of $\overrightarrow V$ with
$\gd\big( x(0)\big) <\gep$, then
$$
 \frac{d}{dt}\Big(\gd\big( x(t)\big)\Big)\geq 
 \sqrt{2\Big(\gep -\gd\big( x(t)\big)\Big)}.
$$
(To be more precise, since $\gd\big( x(t)\big)$ is not necessarily 
differentiable, we should write 
$\liminf_{\gt\to 0}\frac{\gd\big( x(t+\gt)\big)-\gd\big( x(t)\big)}{\gt}$ 
instead of $\frac{d}{dt}\Big(\gd\big( x(t)\big)\Big)$ in the last inequality.)
Thus, for $t=\sqrt{2\Big(\gep -\gd\big( x(0)\big)\Big)}$ we have
$\gd\big( x(t)\big) =\gep$. 

Since $x\in L$ belongs to $\partial L$ iff $\gd (x)=0$, and hence, 
the vector field 
$\overrightarrow V$ points everywhere towards the interior $\textrm{int}(L)$,
it follows from the Peano existence theorem of solution for ordinary 
differential equations, that for any $x\in L$ there is an integral curve 
$c_x(t)$ of $\overrightarrow V$, defined for all $t\geq 0$ with $c_x(0)=x$ and 
with $c_x(t)\in\textrm{int} L$ for $t>0$. Stability of the solutions implies
that $c_x(t)$ depends continuously on $x$ and $t\geq 0$.

As a conclusion, we get that the flow along 
$\overrightarrow V$ for
$\sqrt{2\gep}$ time units defines a deformation retract from $L$ to 
$\mathcal{F}_{\geq\gep}=\{\gd\geq\gep\}.$

If $L$ and $\mathcal{F}$ are
$\gC$-invariant, then $\overrightarrow V(x)$, as it is defined above, is also 
$\gC$-invariant. Hence it induces a vector field on $\gC\backslash L$, 
and a deformation retract from $\gC\backslash L$ to $\gC\backslash\mathcal{F}_{\geq\gep}$.
\end{proof}
 
A second kind of deformation retracts that we shall often use is the following.
Let $A\subset S$ be a closed convex set, and $B\subset S$ a set containing $A$.
We say that $B$ is {\bf star-shaped} 
with respect to $A$ if for any $b\in B$ the geodesic segment 
connecting $b$ to its closest point (the projection) in $A$ is contained 
in $B$. In that case we can define a deformation retract from $B$ to $A$ by 
moving $b$ at a constant rate (equal to the initial distance) 
towards its projection in $A$. We call this the {\bf star-contraction} 
from $B$ to 
$A$. From the convexity of the distance function we conclude that
the star contraction from $B$ to $A$
is distance decreasing, and hence continuous, and if 
the Hausdorff distance $\textrm{Hd}(A,B)$ is finite, then the star-contraction 
gives a homotopy equivalence between $(B,\partial B)$ and $(A,\partial A)$.

More generally, if $A\subset B_1\subset B_2\subset S$, and if the $B_i$'s are closed and
star-shaped with respect to $A$, we can define a deformation retract from $B_2$
to $B_1$ by letting any $b\in B_2\setminus B_1$ flow in the direction of
its projection $P_A(b)$ in $A$ at constant speed $s$, where $s$ equals to
the length of the segment $[b,P_A(b)]\cap (B_2\setminus B_1)$. 
By a similar procedure one can show that there is a deformation retract from $S\setminus B_1$
to $\overline{S\setminus B_2}$.
In this way we obtain:

\begin{lem}[{\bf A generalized star-contraction}]\label{GSC}
Assume that 
\begin{itemize}
\item
$A$ is convex and closed,
\item
$A\subset B_1\subset B_2$ where $B_i$ are closed and star-shaped with respect 
to $A$.
\item
$\textrm{Hd}(A,B_2)<\infty$.
\end{itemize}
Then there is a deformation retract from $(B_2,\partial B_2)$ to $(B_1,\partial B_1)$.
Similarly, there is a deformation retract from $S\setminus B_1$
to $\overline{S\setminus B_2}$.
\end{lem}


\section{Constructing a simplicial complex in a thick submanifold with 
nice boundary}\label{4}

Throughout this section, $M=\gC\backslash S$ is a fixed $S$-manifold with finite volume,
$M'\subset M$ a connected submanifold with boundary, and $\gep >0$ is fixed. 
The main result of this section is Lemma \ref{simplicial-complex}. 

We wish to formulate some convenient conditions on $M'$, under which $M'$ is 
homotopically equivalent to a simplicial complex $\mathcal{R}$ whose 
combinatorics is bounded in terms of $\vol (M)$. 
More precisely, we would like $\mathcal{R}$ to be 
a $\big( d,\ga\cdot\vol (M')\big)$-simplicial complex, where $d$ and $\ga$
are some constants depending only on $S$ and $\gep$ (see Definition \ref{(d,v)-s.c.}).

To construct $\mathcal{R}$ we will use a ``good covering'' argument.
Recall that a cover of a topological space $T$ is called a {\bf good cover} 
if any non empty intersection of sets of the cover is contractible.
In this case the simplicial complex $\mathcal{R}$ which corresponds to
the nerve of the cover is homotopically equivalent to $T$. By definition, the
vertices of $\mathcal{R}$ correspond to the sets of the cover, and a 
collection of vertices form a simplex when the intersection of the 
corresponding sets is non-empty (see \cite{BT}, theorem 13.4). 

In a manifold with injectivity radius bounded uniformly 
from below by $\epsilon$, such a good cover is achieved by taking 
$\epsilon$-balls for which the set of centers form an
$\epsilon /2$ (say) discrete net. In our case, in order to be able to use 
$\gep$-balls (or more precisely $\gep /c$-balls for some constant $c$) 
we shall require that $M'$ lies inside $M_{\geq\gep}$. 
However, $M'$ is not a manifold but a manifold with boundary, and balls may 
not be ``nice'' subsets - they may not be convex or even contractible, and an intersection of
balls may not be connected.
Therefore we should be more careful.
Our problems arise only near the boundary (far away from 
the boundary balls are nice). So we need some control on the geometry of the boundary.

\begin{lem}\label{simplicial-complex}
Let $M=\gC\backslash S$ be a fixed $S$-manifold with finite volume,
let $M'\subset M$ be a connected submanifold with boundary, 
and let $\gep >0$ be fixed. 
\begin{itemize}
\item
Assume that $M'$ is contained in the $\gep$-thick part $M_{\geq \gep}$.
\item
Write $X=M\setminus M'$, and assume that its pre-image $\tilde X\subset S$ 
under the universal covering map is a locally finite union of convex open sets 
with smooth boundary $\ti X=\cup_{i\in I}X_i$
(by locally finite we mean that every compact set in $S$ intersects only finitely many $X_i$'s).
\item
Assume that $M'$ is homotopically equivalent to its $\frac{\gep}{2}$-shrinking 
$)M'(_{\frac{\gep}{2}}$.
\item
Let $b>1$ be a constant which depends only on $S$ and on $\gep$.
\item
Assume that for any point $x\in S\setminus\ti X$ with $d(x,\ti X)\leq\gep$, 
there is a unit tangent vector $\hat n(x)\in T_x(S)$,
such that for each $i\in I$ with $d(x,X_i)=d(x,\ti X)$, 
the inner product of $\hat n(x)$ with the gradient $\nabla D_{X_i}(x)$ satisfies 
$$
 \hat n(x)\cdot\nabla D_{X_i}(x)>\frac{1}{b}.
$$
\end{itemize}
Then there are constants $\ga$ and $d$, depending only on $S$ and on $\gep$, 
such that $M'$ is homotopically equivalent to some 
$\big( d,\ga\cdot\vol (M')\big)$-simplicial complex.
\end{lem}

Throughout this section we use the notation of the statement of 
Lemma \ref{simplicial-complex}.\\

The proof of Lemma \ref{simplicial-complex} relies on a uniform estimate on the distance
between shrinkings of $M'$.

\begin{prop}\label{closeN}
For any $\gt$ and $\gd$ with $\gep\geq \gt+\gd\geq\gt>0$ we have
$$
 \big( )M'(_{\gt+\gd}\big)_{b\gd}\supset )M'(_\gt.
$$
\end{prop}

The proposition can also be stated as follows: For any such $\gt$ and $\gd$, 
the Hausdorff distance between the corresponding sets satisfies
$$
 \textrm{Hd}\big()M'(_{\gt+\gd}~,~)M(_{\gt}\big)= \leq b\gd.
$$
In other words, in order to prove the proposition, we need to show that for any
$x\in )M(_\gt\setminus)M(_{\gt+\gd}=(X)_{\gt+\gd}\setminus (X)_{\gt}$ there is 
$y\notin(X)_{\gt+\gd}$ with 
$d(x,y)\leq b\gd$.

\begin{proof}
Let $x\in(X)_{\gt+\gd}\setminus (X)_{\gt}$. Choose a lifting
$\ti x\in(\ti X)_{\gt+\gd}\setminus (\ti X)_{\gt}$ of $x$.
As $\ti x\in(\ti X)_{\gt+\gd}\setminus (\ti X)_{\gt}$, there is $\gt_1~(\gt+\gd\geq\gt_1\geq\gt)$
such that $\ti x\in\partial(\ti X)_{\gt_1}$.

Let $c(t)$ be the piecewise
geodesic curve of constant speed $b$, passing through $\ti x$, which is defined 
inductively (for $t\geq \gt_1$) as follows:
At time $t_1=\gt_1$ set $c(t_1)=\ti x$ and define its one sided derivative by 
$$
 \frac{d}{dt}_+ c(t_1)=b\hat n(\ti x),
$$ 
and define 
$$
 I_c(t_1)=\{ i\in I:x\in\partial (X_i)_{t_1}\}.
$$
Identify $c$ with the constant speed geodesics determined by this condition, as
long as $t >t_1$ and $c(t )\notin \overline{(X_i)_\gt}$ for any $i\notin 
I_c(t_1)$. Let $t_2$ be the first time (if such exist) $t_1\leq t<\gt+\gd$ at 
which $c(t)$ hits some $\overline{(X_i)_t}$ for $i\notin I_c(t_1)$. 
As the collection $\{ X_i\}_{i\in I}$ is locally finite, $t_2$ is well defined and strictly 
bigger than $t_1$. We claim that $c(t_2)\in\partial (\tilde X)_{t_2}$, and 
that $c(t)\notin (\ti X)_t$ for $t_1<t<t_2$. To prove this, we
need to show that $c(t )\notin (\ti X_i)_t$ for any 
$t_1\leq t\leq t_2$ and $i\in I_c(t_1)$. Fix such $i$ and observe that
$$
 \frac{d}{dt}|_{t=t_1} D_{X_i}\big(c(t)\big)=
 b\hat n(x)\cdot \nabla D_{X_i}(x)>b\frac{1}{b}=1.
$$ 
Since the convex function $D_{X_i}\big( c(t )\big)$ has 
non-decreasing derivative we get that $\frac{d}{dt}\Big( D_{X_i}\big( c(t)\big)\Big) >1$ 
for any $t_2>t>t_1$, and hence, 
$D_{X_i}\big( c(t)\big)>D_{X_i}\big( c(t_1)\big)+(t-t_1)=t$ 
for any such $t$. 
Thus, the point $c(t)$ is outside $\big( (X_i)_{t_1}\big)_{t -t_1}=(X_i)_t$.

Then we define the second piece of $c$ by the condition that its one sided
derivative $\frac{d}{dt}_+ c(t_2)$ satisfies 
$\frac{d}{dt}_+ c(t_2)=b\hat n\big( c(t_2)\big)$, 
and we define 
$$ 
I_c(t_2)=\{ i\in I:c(t_2)\in \partial (X_i)_{t_2}\}.
$$

Note that since $t_2$ is smaller than $\gt+\gd <\gep$ and 
$c(t_2)\in\partial (X_2)$, the  direction $\hat n\big( c(t_2)\big)$ is well 
defined. In this way, we continue to define $c$ inductively.

If $t_i$ converge to some $t_{\omega_0}<\gt+\gd$, then we define 
$c(t_{\omega_0})$ to be the limit of $c(t_i)$ and 
$\frac{d}{dt_+}c(t_{\omega_0})=b\hat n\big( c(t_{\omega_0})\big)$
and $I_c(t_{\omega_0})=\{i\in I:c(t_{\omega_0})\in\partial 
(X_i)_{t_{\omega_0}}\}$,
and denote the next index by $\omega_0+1$, and so on.
By this way we obtain a piecewise
geodesic curve with at most countably many pieces, connecting 
$\ti x=c(\gt_1)$ to $\ti y=c(\gt +\gd)$, 
with $\ti y=c(\gt +\gd )\notin (\tilde X)_{\gt +\gd}$. 
Since the length of 
$c\big( [\gt_1,\gt +\gd ]\big)$ is $b(\gt +\gd-\gt_1)\leq b\gd$ 
this proves the proposition.
\end{proof}

\begin{cor}\label{cover}
Let $\gt$ and $\gd$ be as in Proposition \ref{closeN}. Let $\mathcal{C}$ be a
collection of balls of radius $(b+1)\gd$ for which the set of centers
$\mathcal{C}'$ form a maximal $\gd$ discrete subset of $)M(_{\gt+\gd}$.
Then $\mathcal{C}$ is a cover of $)M(_\gt$.
\end{cor}

\begin{proof}
As $\mathcal{C}'$ is maximal $(\mathcal{C}')_\gd\supset )M(_{\gt+\gd}$.
By Proposition \ref{closeN} we have 
$$
 \big( )M'(_{\gt+\gd}\big)_{b\gd}\supset )M'(_\gt.
$$
Thus
$$ 
 \cup_{C\in \mathcal{C}}C=(\mathcal{C}')_{(b+1)\gd}=
 \big( (\mathcal{C}')_\gd\big)_{b\gd}\supset
 \big( )M'(_{(\gt +\gd)}\big)_{b\gd}\supset )M'(_\gt.
$$
\end{proof}

In the sequel we will take $\gt=\frac{\gep}{2}$.
The next proposition provides a uniform bound on the curvature of the smooth pieces
of $\partial )M'(_{\frac{\gep}{2}}$.

\begin{prop}\label{smallC}
For any $X_i$ and any point $x \in \partial (X_i)_{\frac{\epsilon}{2}}$,
the $\frac{\gep}{2}$-ball, whose boundary sphere is tangent at $x$ to  
$\partial (X_i)_{\frac{\epsilon}{2}}$, which lies on the same side of 
$\partial (X_i)_{\frac{\gep}{2}}$ as $(X_i)_{\frac{\gep}{2}}$,
is contained in $(X_i)_{\frac{\epsilon}{2}}$.
\end{prop}

\begin{proof}
The distance between $x$ and its closest point $p_{X_i}(x)$ in the closed 
convex set $X_i$ is easily seen to be $\frac{\gep}{2}$, and the 
$\frac{\gep}{2}$-ball 
centered at $p_{X_i}(x)$ is the required one.
\end{proof}

The following proposition follows directly from the definition of a deformation 
retract.

\begin{prop}\label{B-B'}
Let $B'\subset B$ be topological spaces, and let $F_t~\big( t\in [0,1]\big)$ 
be a deformation retract of $B$ such that $F_t(b)\in B'$ for any 
$b\in B',t\in [0,1]$. 
Then $F_t|_{B'}$ is a deformation retract of $B'$.
\end{prop}

Let $B_r$ denote an arbitrary ball in $S$ of radius $r$, and $B_r(x)$ the ball 
of radius $r$ centered at $x$. The following proposition follows from the fact 
that the volume of a ball of radius $r$ is a convex function of $r$ (because 
the surface area of the $r$-sphere is an increasing function of $r$).

\begin{prop}
There is a constant $m$ such that for any $\delta <1$, 
$$
 m\cdot \vol (B_{\delta /2})
 \geq \vol (B_{(b+1.5)\delta}).
$$ 
Thus any $\delta$-discrete subset of 
$B_{(b+1)\delta}$ consists of at most $m$ elements.
\end{prop}

For a finite set
$\{ y_1,...,y_t\} \subset S$ we denote by $\sigma (y_1,...,y_t)$
its Chebyshev center, i.e. the unique point $x$ which minimizes the function
$\max_{1\leq i\leq t}d(x,y_i)$.

We will soon take $B$ to be an intersection of balls, and $B'\subset B$ to be
the intersection of $B$ with $)M'(_{\frac{\gep}{2}}$. 
We intend to use \ref{B-B'} in order to show that under some certain 
conditions, $B'$ is contractible.
It will be natural to use the star contraction to the Chebyshev center 
of the centers of the associated balls - The deformation retract which 
makes any point of $B$ flow along the geodesic segment which connects it to 
the required Chebyshev center, at constant velocity $s$, where $s$ 
equals the initial distance. In order to do this, we need the following:

\begin{prop}[Defining the constant $\gd$] \label{positiveD}
There exists $\gd~(0<\delta < \frac{\epsilon}{2(b+1)})$ 
such that for any point $x\in S$, any $\frac{\gep}{2}$-ball
$C$ which contains $x$ on its boundary sphere, and  
any $m$ points $y_1,\ldots ,y_m\in B_{(b+1)\delta }(x)\setminus 
(C)_{\delta}$, the inner product of the external normal vector of $C$ at $x$ 
with the tangent at $x$ to the geodesic segment $[x,\sigma (y_1,\ldots ,y_m)]$ 
is positive.
\end{prop}

\begin{proof}
Assume the contrary. Then there is a sequence $\delta_n \to 0$,
a corresponding sequence $(C_n)$ of $\frac{\epsilon}{2}$-balls tangent to some 
$x_0\in S$ (which we may assume converge to some fixed such ball), and a 
corresponding sequence of
$m$-tuples of points 
$y_1^n,\ldots ,y_m^n\in B_{(b+1)\delta_n}(x_0)\setminus ( C_n)_{\gd_n}$, 
such that the inner product of the
tangent at $x_0$ to the geodesic segment $[x_0,\sigma (y_1^n,...,y_m^n)]$ with
the external normal vector of the corresponding sphere $\partial C_n$ is non-positive (we may fix $x=x_0$ 
since $G$ acts transitively).

Rescaling the metric each time we can
assume that $\delta_n$ is fixed and equals $1$. We then get a sequence of 
Riemannian metrics converging, on the ball of radius $b+1$ around $x_0$,
to the Euclidean metric on the ball of radius $b+1$. 
More precisely, we look at the ball of radius $b+1$ in the tangent space
$T_{x_0}S$. We identify it
each time with the ball of radius $(b+1)\delta_n$ around $x_0$ in $S$ via the 
map $X\to \exp_{x_0}(\delta_n X)$, and we rescale the metric there to 
$$
 d_n(X,Y)=\frac{1}{\delta_n}\cdot d\big(\exp_{x_0}(\delta_nX),\exp_{x_0}(\delta_nY)\big).
$$
( All these metrics induce the same topology.)

Now, in the rescaled metrics our tangent balls tend to a half space 
(since $\epsilon /\delta_n \to \infty$), and we may as well assume that our 
$m$-tuples also converge. 
In the limit, we get an $m$-tuple of points
in a Euclidean space at distance at least $1$ (and at most $b+1$)
from a half space, for which
the inner product of the external normal vector to this half space with the 
vector $\overrightarrow v$, pointing from \textrm{so}me $x_0$ on the boundary 
hyper-plane, towards the Chebyshev center of this $m$-tuple is non-positive. This is an absurd.
\end{proof}  

Finally, we claim

\begin{prop}
Let $\mathcal{C}$ be a collection of balls of radius $(b+1)\delta$, for which 
the set of centers $\mathcal{C}'$ form a maximal $\delta$-discrete subset of
$)M'(_{\frac{\gep}{2} +\gd}$. Then $\mathcal{C}$, i.e. the restrictions of
its sets to $)M'(_{\frac{\gep}{2}}$, is a good cover of 
$)M'(_{\frac{\gep}{2}}$.
\end{prop}

\begin{proof}
By corollary \ref{cover}, $\mathcal{C}$ is a cover of $)M'(_{\frac{\gep}{2}}$.

Let $B$ be the intersection of $m$ (not necessarily different) balls of 
$\mathcal{C}$, with centers $y_1,\ldots,y_m\in\mathcal{C}$. 

Proposition \ref{smallC} implies that for any 
$\tilde x\in \partial )\tilde M'(_{\frac{\gep}{2}}$ 
and for any  $X_i$ with
$\tilde x\in \partial  (X_i)_{\frac{\epsilon}{2}}$
there is an $\frac{\epsilon}{2}$-ball,
tangent to $\partial (X_i)_{\frac{\gep}{2}}$ at $\tilde x$,
which is contained in $(X_i)_{\frac{\gep}{2}} \subset 
(\tilde X)_{\frac{\gep}{2}}$.
Thus, if in addition the image $x$ of $\tilde x$ belongs to $B$, then 
Proposition \ref{positiveD} implies that the geodesic segment 
$[x,\sigma (y_1,\ldots ,y_m)]$ lies inside $B'=B\cap )M'(_{\frac{\gep}{2}}$.
Let us explain this point as follows. Let $c(t),~t\in [0,1]$ be a 
parameterization of the geodesic segment $[x,\sigma (y_1,\ldots ,y_m)]$. 
Propositions \ref{smallC} and \ref{positiveD} imply that 
$c(t)\notin \overline{(X)}_{\frac{\gep}{2}}$ for all sufficiently small $t>0$.
Assume that $c(t)\in\partial (X)_{\frac{\gep}{2}}$ for some $t_0$, $0<t_0<1$, 
(and let $t_0$ be the first such  time). Then Propositions \ref{smallC} and 
\ref{positiveD} applied to $c(t_0)$ imply that 
$c(t_0-\gD t)\in(X)_{\frac{\gep}{2}}$ for small $\gD t$. But this is a 
contradiction.
We conclude that if $B$ is not empty 
then $\sigma (y_1,\ldots ,y_m)\in B'$ and the star-contraction from $B$ to
$\sigma (y_1,\ldots ,y_m)$ induces a contraction of $B'$. Hence the 
set $B'$ is non-empty and contractible.
This means that $\{ C':C\in \mathcal{C} \}$ is a good cover of 
$)M'(_{\frac{\gep}{2}}$, where $C':=C\cap )M'(_{\frac{\gep}{2}}$.
\end{proof}

We conclude that $)M'(_{\frac{\gep}{2}}$, and therefore $M'$, is homotopically 
equivalent to the simplicial complex $\mathcal{R}$ which corresponds to the 
nerve of $\mathcal{C}$. Since the collection of centers $\mathcal{C}'$ is 
$\gd$-discrete, and therefore the $\frac{\gd}{2}$-balls with the same centers
are disjoint, we get that $|\mathcal{C}'|$,  
the number of the vertices of $\mathcal{R}$, is 
$\leq\frac{\vol (M)}{\vol (B_{\gd /2})}$. 
Since the sets of $\mathcal{C}$ are subsets of $(b+1)\gd$-balls, each of them 
intersects at most 
$\frac{\vol ( B_{2(b+1)\gd+\gd /2})}{\vol (B_{\gd /2})}$ of the 
others. Hence, the degree of any vertex in $\mathcal{R}$ is 
$\leq d:=\frac{\vol ( B_{2(b+1.25)\gd})}{\vol (B_{\gd /2})}$.
This completes the proof of Lemma \ref{simplicial-complex}.


\section{An arithmetic variant of the Margulis lemma}\label{AML}

The classical Margulis lemma yields information on the algebraic structure
of a discrete group of isometries which is generated by ``small elements''. If,
in addition, this discrete group lies inside an arithmetic group then we can 
say a little more. In this section we shall explain this idea in the
non-uniform case. 

The Lie group of isometries $G=\textrm{Isom}(S)$
is center-free, semi-simple, without compact factors and with 
finitely many connected components. Its identity component $G^0$ coincides 
with the connected component of the group of real points $\Bbb{G}(\BR )^0$ of some
$\BQ$-algebraic group $\Bbb{G}$. 
We shall identify $G^0$ with its adjoint group
$\text{Ad}(G^0)\leq\text{GL}(\mathfrak{g})$ and think of it as a group of 
matrices. We will denote by $\mu$ a fixed Haar measure on $G$.

\begin{lem}[\bf An arithmetic variant of the Margulis lemma]\label{arith-M-L}
There are constants $\gep =\gep (S)>0$ and $m=m(S)\in\BN$, such that 
if $\gC\leq G$ is a non-uniform torsion-free arithmetic lattice, then 
for any $x\in S$, the group of real points of the Zariski closure
$\big(\overline{\gC_\gep (x)}^z\big)_\BR$ of the group 
$$
 \gC_\gep (x)=\langle\gc\in\gC:d_\gc (x)\leq\gep\rangle
$$
has at most $m$ connected components, and its identity component is a 
unipotent subgroup.
\end{lem}

The following claims (\ref{Q-structure}, \ref{5.3}, \ref{5.4}) 
are well known (c.f. \cite{Mar1}, chapter $IX$ section 4). 

\begin{lem}\label{Q-structure}\label{5.2}
For any non-uniform arithmetic lattice $\gD\leq G^0$, there is a rational 
structure on $G^0$ (coming from a rational structure on the vector space $\mathfrak {g}$)
with respect to which $\gD$ is contained in $G^0(\BQ )$ and
commensurable to $G^0(\BZ )$. Conjugating by an element $g\in G^0(\BQ )$ we can
assume that $\gD\subset\ G^0(\BZ )$.
\end{lem}

\begin{proof}[Explanation]
In general, if $\gD$ is an arithmetic lattice in $G$, there is 
a compact extension $H$ of $G$, a rational structure on $H$, and a subgroup
$\gD'\leq H_\BQ$ commensurable to $H_\BZ$ whose projection to $G$ is $\gD$.   
However, when $\gD$ is non-uniform we can always take $H=G$.
This could be deduced, for example, from the fact that $\gD$ is almost 
generated by unipotent elements, and that compact groups contain no non-trivial 
unipotent element.

For a given rational structure, the group $G_\BZ$ is defined only up to 
commensurability. A subgroup $\gD\leq G_\BQ$ which is commensurable to 
$G_\BZ$ is always conjugate, by an element of $G_\BQ$, to a subgroup of $G_\BZ$.
\end{proof}

If $g\in G$ is close to $1\in G$, then its characteristic polynomial (i.e. the characteristic 
polynomial of the endomorphism $\text{Ad}(g)$) is close to $(\gl -1)^n$ where 
$n=\dim (\mathfrak {g})$. In particular:

\begin{prop}\label{5.3}
There is an identity neighborhood $\gO_1\subset G$ such that, if $g\in \gO_1$,
and the characteristic polynomial of $g$ has integral coefficients, 
then $g$ is unipotent.
\end{prop}

\ref{5.2} and \ref{5.3} implies:

\begin{cor}[\cite{Mar1} 4.21]\label{5.4}
For any non-uniform arithmetic lattice $\gD\leq G$, the intersection 
$\gO_1\cap\gD$ consists of unipotent elements only.
\end{cor}

Recall also the following theorem of Zassenhaus and Kazhdan-Margulis (see \cite{Rag} 
theorem 8.16):   

\begin{thm}[Zassenhaus, Kazhdan-Margulis]
There exists an identity neighborhood $\gO_2\subset G$, called a
Zassenhaus neighborhood, so that for any discrete subgroup $\gS\leq G$,
the intersection $\gS\cap\gO_2$ is contained in a connected nilpotent Lie
subgroup of $G$.
\end{thm}

\begin{proof}[Proof of Lemma \ref{arith-M-L}]
Let $\gO\subset G$ be a relatively compact symmetric identity neighborhood 
which satisfies $\gO^2\subset\gO_1\cap\gO_2$.

Fix an integer $m$,
$$
 m> \inf_{h\in G}\frac{\mu \big(\{g\in G:d_g(x)\leq 1\}\cdot h\gO h^{-1} \big)}
 {\mu (\gO )}, 
$$ 
and
$$
 \gep =\frac{1}{m}.
$$ 
As $G$ is unimodular, $m$ can be chosen independently of $x$. 
Replacing $\gO$ by some conjugate $h\gO h^{-1}$, if needed, we can assume that
$$
 m>\frac{\mu \big(\{g\in G:d_g(x)\leq 1\}\cdot\gO \big)}{\mu (\gO )}.
$$ 
Let
$$
 \gC_\gO=\langle\gC\cap\gO^2\rangle.
$$
Then 
$$
 [\gC_\gep (x): \gC_\gep (x)\cap\gC_\gO ]\leq m.
$$
To see this, assume for a moment that this index was $\geq m+1$. Then we 
could find $m+1$ representatives $\gc_1,\gc_2,\ldots,\gc_{m+1}\in\gC_\gep (x)$
for different cosets of $\gC_\gep (x)\cap\gC_\gO$ in the ball of radius $m$
in $\gC_\gep (x)$ according to the word metric with respect to the generating
set $\{\gc\in\gC_\gep (x): d_\gc (x)\leq\gep\}$. As they belong to different
cosets, $\gc_i\gO\cap\gc_j\gO =\emptyset$ for any $1\leq i<j\leq m+1$.
Since $d_{\gc_i}(x)\leq m\cdot\gep =1$ these $\gc_i$'s
are all inside $\{ g\in G:d_g(x)\leq 1\}$. This contradicts the assumption 
$m\cdot\mu (\gO )>\mu \big(\{ g\in G:d_g(x)\leq 1\}\cdot\gO\big)$.

It follows from the Zassenhaus-Kazhdan-Margulis theorem that $\gC_\gO$ is contained in a 
connected nilpotent Lie subgroup of $G$, and therefore, by Lie's theorem,
$\gC_\gO$ is triangulable over $\BC$. As $\gC_\gO$ is generated by unipotent
elements, it follows that $\gC_\gO$ is a group of unipotent elements.
Thus the Zariski closure $\overline{\gC_\gO}^z$ is a unipotent algebraic group,
and hence, the group of its real points $(\overline{\gC_\gO}^z)_\BR$ is 
connected in the real topology. Similarly, its subgroup
$(\overline{\gC_\gep (x)\cap\gC_\gO}^z)_\BR$ is connected.
Clearly, $(\overline{\gC_\gep (x)\cap\gC_\gO}^z)_\BR$ is the identity 
connected component of $(\overline{\gC_\gep (x)}^z)_\BR$, 
and its index is at most $m$.

\end{proof}

\begin{rem}
Although $\gep$ could be taken to be $1/m$, we use different letters for them
because they play different roles.
In the sequel we will assume that the above lemma is satisfied with $\gep$
replaced by $10\gep$.
\end{rem}

\begin{rem}\label{5.7}
It follows from Lemma \ref{arith-M-L} that if $\gc$
is an element of a non-uniform arithmetic lattice of $G$ and 
$\inf d_\gc < \gep$, then $\gc^j$ is unipotent for some $j\leq m$.
This implies that in a non-compact arithmetic $S$-manifold there are no closed geodesic of length 
$\leq\gep$, and that in the $\gep$-thick thin decomposition, the thin part has no
compact connected component. For example, the $\gep$-thin part of any 
non-compact arithmetic hyperbolic surface (or more generally of any rank-1
manifold) consists only of cusps.
\end{rem}


\section{The proof of \ref{thmA}(1)}\label{6}

In this section we shall prove:

\begin{thm}[{\bf Theorem \ref{thmA}(1) of the introduction}]\label{thm1}
Let $S$ be a symmetric space of non-compact type. Then there are constants
$\ga$ and $d$ (depending only on $S$) such that any non-compact arithmetic
$S$-manifold $M$ is homotopically equivalent some 
$\big( d,\ga\cdot\vol (M)\big)$-simplicial complex.
\end{thm}

Fix $\gep' =\gep' (S),m=m(S)$, such that Lemma \ref{arith-M-L} is satisfied 
with $10\gep'$, and $m$. Assume that $M=\gC\backslash S$ is a given non-compact 
arithmetic $S$-manifold. Denote by $\gC^u$ the set of unipotent elements in 
$\gC$,
$$
 \gC^u=\{ \gc\in\gC :\gc \textrm{ is unipotent}\},
$$
and define
$$
 \tilde X=\cup_{\gc\in\gC^u\setminus\{ 1\}}\{ d_\gc<\gep'\},
$$
and
$$
 \tilde M'=\tilde M \setminus \tilde X,
$$
and let $X\subset M$ and $M'\subset M$ be the images of $\ti X$ and 
$\ti M'$ under the universal covering map.

In order to prove Theorem \ref{thm1} we shall show:

\begin{enumerate}
\item\label{1}
There is a deformation retract from $M$ to $M'$. In particular $M$ is 
homotopically equivalent to $M'$.
\item\label{2}
$M'$ satisfies the conditions of Lemma \ref{simplicial-complex} with 
$\gep =\gep'/m$.
\end{enumerate}

Let us start with some preliminaries. Let $W\subset S(\infty )$ be a Weyl
chamber of the Tits spherical building. Then $W$ is isometric to a spherical
simplex.
We define its center of mass $z$ to be
$$
 z=\frac{\int_W xd\mu (x)}{\|\int_W xd\mu (x)\|}
$$
where $\mu$ is the Lebesgue measure on the sphere. Since $W$ is 
contained in a half 
sphere (see \cite{BGS} appendix 3), we get $\|\int_W xd\mu(x)\|\neq 0$. Since $W$
is convex (in the spherical metric), and since its interior is non-empty and is 
exactly the set of regular points in $W$, it follows that $z$ is regular
and contained in $W$.

We shall also use:

\begin{lem}\label{321}
Assume that $g$ is a parabolic isometry (i.e. $\min (g)=\emptyset$), $c(t)$ is a regular geodesic, and
$g$ stabilizes $c(-\infty )$. Then $\frac{d}{dt}d_g\big( c(t)\big)\neq0$
for any $t\in\BR$.
\end{lem}

\begin{proof}
Let $c(t)$ be a regular geodesic in $S$, and let $g$ be an isometry which 
stabilizes $c(-\infty )$. Assume that 
$\frac{d}{dt}|_{t=t_o}d_g\big( c(t)\big) =0$, then the function 
$d_g\big( c(t)\big)$, being analytic and convex 
must be constant, since it doesn't tend to $\infty$ as 
$t\to -\infty$. This means that $g\cdot c$ is a geodesic parallel to $c$.
We conclude that $g$ preserves the unique flat which contains $c$, since this
flat is exactly the set of points through which there is a geodesic parallel
to $c$. But this implies that $g$ is semisimple, a contradiction.
\end{proof}

Furthermore, if $g$ is unipotent, then $d_g\big( c(t)\big)\to 0$ as $t\to -\infty$
(see \cite{BGS} appendix 5, section 4).

\begin{proof}[{\bf The proof of (\ref{1})}]
We will show that the conditions of Lemma \ref{deformation-retract} are 
satisfied with $L=Y=S$ and the $\gC$-invariant family of functions
$\mathcal{F}=\{d_\gc\}_{\gc\in\gC^u\setminus\{ 1\}}$. 
Then $\ti M'=\mathcal{F}_{\geq\gep'}$.

The finiteness of the sets 
$$
 \Psi_{x,\gt }=\{\gc\in\gC^u\setminus\{ 1\}:d_\gc (x)\leq\gt\}
$$ 
follows from the compactness of $\{ g\in G:d_g(x)\leq\gt\}$
together with the discreteness of $\gC$. All the functions 
$\{ d_\gc\}_{\gc\in\gC^u\setminus\{ 1\}}$ are convex, strictly positive, 
without critical values. 
We need to define the continuous function $\gb :(0,3\gep )\to\BR^{>0}$, 
and the appropriate direction $\hat n(x)\in T_x(S)$ for any $x\in S$ with
$\Psi_{x,3\gep}\neq\emptyset$.

Fix $x\in S$ with $\Psi_{x,3\gep}\neq\emptyset$ and let 
$\Psi_{x,3\gep}=\{\gc_1,\gc_2,\ldots,\gc_k\}$.
By Lemma \ref{arith-M-L} the Zariski closure $\overline\gD^z$ of
the group $\gD =\langle \gc_1,\gc_2,\ldots,\gc_k \rangle$ 
has a unipotent identity component, and $\leq m$ connected components. 
Since $\gc_i$ is unipotent,
it is contained in the Zariski closure of the cyclic group generated by any 
power of it. 
As $\gc_i^j$ belongs to the identity component $(\overline{\gD}^z)^0$ for some
$j\leq m$, also $\gc_i\in(\overline{\gD}^z)^0$. Since this holds for any
generator $\gc_i$, we get that
$\overline{\gD}^z=(\overline{\gD}^z)^0$ and hence the group of its
real points $(\overline{\gD}^z)_\BR$ is a connected unipotent group.
Hence $\gc_1,\gc_2,\ldots,\gc_k$ are contained in a connected unipotent group.

Let $N\leq G$ be a maximal 
connected unipotent subgroup which contains $\gc_1,\gc_2,\ldots ,\gc_k$.
Let $W\leq S(\infty )$ be the Weyl chamber of the Tits boundary of $S$ 
which corresponds to $N$. Let $c(t)=c_x(t)$ be the
geodesic line with $c(0)=x$ for which $c(-\infty )$ is the center of mass of
$W$.

By Lemma \ref{321}
$$
 \frac{d}{dt}|_{t=0}\Big( d_g\big( c(t)\big)\Big)>0
$$ 
for any $g\in N\setminus\{ 1\}$. In addition, the continuous function
$$
 h(g)=\dot c(0)\cdot\nabla d_g(x)
$$ 
attains a minimum on the compact set
$$
 \{g\in N:d_g(x)=\gt\}.
$$
We define $\gb (\gt )$ to be this minimum (then $\gb (\gt )$ is defined for
any $\gt >0$). By definition, 
$\dot c(0)\cdot\nabla d_{\gc_i}(x)\geq\gb\big(d_{\gc_i}(x)\big)$.  
Moreover, it is easy to see that $\gb$ is a 
continuous positive function. Since $G$ acts transitively on the
set of couples $(W,x)$ of a Weyl chamber $W\subset S(\infty )$ and a point $x\in S$,
$\gb$ is independent of $N$ and of $x$.
The conditions of Lemma \ref{deformation-retract} are satisfied with the 
tangent vector $\hat n(x)=\dot c_x(0)$.
\end{proof}

\begin{rem}
We used the existence and uniqueness of the center of mass
of $W$ in order to define 
$\gb$ in a canonical way. However, as any two Weyl chambers are isometric,
we could choose arbitrarily a regular point in one Weyl chambers, and translate
it to any other Weyl chambers, and to use these points when defining 
$c(-\infty )$ and then $c, \hat n(x)$ and
$\gb$ in our proof.
\end{rem}

\begin{rem}\label{finite-presentability}
A slight modification of the argument above, yields an elementary proof of the 
fact that $M=\gC\backslash S$ is homotopic to a compact manifold with boundary, and 
hence that $\gC$ is finitely presented, when $\gC\leq G$ is any torsion free (non-uniform)
lattice. We used the arithmeticity of $\gC$ in order to get a uniform estimate for all $\gC$'s.
However, for a fixed $\gC$, we could have used corollary 11.18 from \cite{Rag} instead of 
corollary \ref{5.4} above, and by this to avoid the assumption that $\gC$ is arithmetic. Moreover, 
we don't really have to assume that $\gC$ is torsion free (when proving just the finiteness without
explicit estimate). The same argument shows that in general, any $S$-orbifold 
$\gC\backslash S$ is homotopic to a compact orbifold with boundary, and it is not hard to show by the same 
means that $\gC\backslash G$ is homotopic to a compact manifold with boundary. In particular, our method 
gives a quite elementary proof of the finite presentability of lattices.
Furthermore, changing appropriately the factor $\sqrt{\cdot}$ in the vector field which induces
the deformation retract, so that it will decay more slowly, we can get a simple proof that 
$\gC\backslash G$ is diffeomorphic to the interior of a compact manifold with boundary.
\end{rem}

\begin{proof}[{\bf The proof of (\ref{2})}]
In order to check that the conditions of Lemma \ref{simplicial-complex}
are satisfied, we will show that

\begin{itemize}
\item $M'\subset M_{\geq\gep}$ (where $\gep =\frac{\gep '}{m}$).
\item $\ti X$ is a locally finite union of convex open sets with smooth boundary.
\item $M'$ is homotopically equivalent to its $\frac{\gep'}{2m}$-shrinking 
$)M'(_{\frac{\gep'}{2m}}$.
\end{itemize}
And that there are
\begin{itemize}
\item $b=b(S)>0$ and
\item $\hat n(x)\in T_x(S)$ for any 
$x\in (\ti X)_{\frac{\gep'}{2m}}\setminus (\ti X)$
\end{itemize}
such that
$$
\hat n(x)\cdot\nabla D_{\{\gc\leq\gep'\}}(x)>\frac{1}{b}
$$
for any $\gc\in\gC^u\setminus\{ 1\}$ with 
$D_{\{\gc\leq\gep'\}}(x)\leq\frac{\gep'}{2m}$ (this is stronger than the condition on the
inner products required in Lemma \ref{simplicial-complex}).\\

If $x\in S$ and $d_\gc (x)\leq\frac{\gep'}{m}\leq\gep'$ for some $\gc\in\gC$, 
then by Lemma \ref{arith-M-L} $\gc^j$ is unipotent for some $j\leq m$,
and since $d_{\gc^j} (x)\leq\frac{j\cdot\gep'}{m}\leq\gep'$, we obtain that
$x\in\ti X$. This proves that $M'$ is contained in the $\frac{\gep'}{m}$-thick
part $M_{\geq\frac{\gep'}{m}}$.\\

Since for $\gc\in\gC^u\setminus\{ 1\}$ the displacement function is convex and
analytic and $\inf d_\gc=0$, and since $\gC$ is discrete, we have that 
$\ti X=\cup_{\gc\in\gC^u\setminus\{ 1\}}\{d_\gc<\gep\}$ is a locally finite 
union of convex open sets with smooth boundary.
\\

We shall prove that there is a deformation retract from $M'$ to 
its $\frac{\gep'}{2m}$-shrinking $)M'(_{\frac{\gep'}{2m}}$
by showing that the conditions of Lemma \ref{deformation-retract}
are satisfied with 
${\mathcal F}=\{ D_{\{ d_\gc\leq\gep'\}}\}_{\gc\in\gC^u\setminus\{ 1\}}$,
$Y=S$, $L={\mathcal F}_{\geq 0}=\ti M'$, and 
${\mathcal F}_{\geq\frac{\gep'}{2m}}=)\ti M'(_{\geq\frac{\gep'}{2m}}$.
(In this proof $\frac{\gep'}{2m}$ plays the rule of $\gep$ in Lemma
\ref{deformation-retract}.)

The finiteness of the sets 
$$
 \Psi_{x,\gt}=\{\gc\in\gC^u\setminus\{ 1\}:
 D_{\{d_\gc\leq\gep'\}}(x)\leq\gt\}
$$ 
follows from the discreteness of $\gC$
together with the compactness of 
$$
\{g\in G:D_{\{d_g\leq\gep\}}(x)\leq\gt\}.
$$ 

We shall define the direction $\hat n(x)\in T_x(S)$ analogously 
to the way it is done in the proof of (\ref{1}), 
and shall show that the conditions of Lemma \ref{deformation-retract} are 
satisfied with the constant function $\frac{\gb (\gep' )}{2}$, where $\gb$ 
is the function defined in the proof of (\ref{1}).

Fix $x\in L=\mathcal{F}_{\geq 0}$, with $\Psi_{x,\frac{3\gep'}{2m}}\neq\emptyset$, and denote
$\Psi_{x,\frac{3\gep'}{2m}}=\{\gc_1,\gc_2,\ldots,\gc_k\}$.
Since $\gep'+2\frac{3\gep'}{2m}=4\gep'$ we have
$(\{ d_{\gc_i}\leq\gep'\} )_{\frac{3\gep'}{2m}}\subset\{d_{\gc_i}\leq 4\gep'\}$
which implies that $\cap_1^k\{d_{\gc_i}\leq 4\gep'\}\neq\emptyset$ and hence
$\gc_1,\gc_2,\ldots,\gc_k$ are contained in a connected 
unipotent group. As in the proof of (\ref{1}), let $N$ be a 
maximal connected unipotent group which contains $\gc_1,\gc_2,\ldots,\gc_k$, 
let $W\subset S(\infty )$ be the Weyl chamber which corresponds to $N$ in the 
Tits boundary. Let $c(t)=c_x(t)$ be the geodesic line with $c(0)=x$ and $c(-\infty)=$ the center 
of mass of $W$. Since $\gc_i\in N$,
$d\big( c(t),\gc_i\cdot c(t)\big)$ tends to $0$ as $t\to -\infty$, and since
$d_{\gc_i}(x)>\gep'$ we have
$d_{\gc_i}\big( c(t_0)\big)=\gep'$ for some negative $t_0$. 

The function $d_{\gc_i}\big( c(t)\big)$ is convex and hence has a non-decreasing derivative,
thus, taking $\hat n(x)=\dot c(0)$, we have:

\begin{eqnarray*}
  \hat n(x)\cdot\nabla D_{\{d_{\gc_i}\leq\gep'\}}
 (x)
&=& \dot c(0)\cdot\nabla D_{\{d_{\gc_i}\leq\gep'\}}(x)\\ 
& =& 
 \dot c(t_0)\cdot\frac{\nabla d_{\gc_i}\big( c(0)\big)}
 {\|\nabla d_{\gc_i}\big( c(0)\big)\|}\\
&\geq &
 \dot c(0)\cdot\frac{\nabla d_{\gc_i}\big( c(0)\big)}{2}\\
&=& 
\frac{1}{2}\cdot\frac{d}{dt}|_{t=0}d_{\gc_i}\big( c(t)\big)\\
& \geq& 
\frac{1}{2}\cdot\frac{d}{dt}|_{t=t_0}d_{\gc_i}\big( c(t)\big)\\
&=&
\frac{1}{2}\cdot\dot c(t_0)\cdot\nabla d_{\gc_i}\big( c(t_0)\big)
 \geq\frac{\gb(\gep')}{2},
\end{eqnarray*}
where $\gb$ is the function defined in the proof of 
\ref{1}.
In the above computation we made use of the facts that 
$\|\nabla  D_{\{d_{\gc_i}\leq\gep'\}}\| =1$ everywhere outside 
$\{d_{\gc_i}<\gep '\}$, 
and that $d_{\gc_i}$ is 2-Lipschitz and hence $\|\nabla d_{\gc_i}\|\leq 2$.

We completed the verification that the conditions of Lemma 
\ref{deformation-retract} are satisfied, and hence proved that $M'$ is 
homotopic to its $\frac{\gep'}{2m}$-shrinking.\\

Finally, observe that we can take also the constant $b$
to be $b=\frac{2}{\gb(\gep')}$, and the unit tangent vector $\hat n(x)$ from 
\ref{simplicial-complex} to be the same $\hat n(x)$ used above. Then all the conditions of Lemma
\ref{simplicial-complex} are satisfied, and the proof of (\ref{2}) is 
completed. 
\end{proof}


\section{Estimating angles at corners of the boundary}\label{bigangles}\label{7}

The next 3 sections are devoted to the remaining cases (2), (3) of Theorem \ref{thmA}.
The main result on this section is Theorem \ref{nice-boundary}.

Our method is to replace a manifold $M$ by a submanifold with boundary $M'$, 
in which the injectivity radius is uniformly bounded from below, and then  
to apply Lemma \ref{simplicial-complex}. This requires some
information on the boundary. The pre-image $\ti X$ of the complement 
$X=M\setminus M'$ in the universal covering is required to be
a locally finite union of 
convex open sets with smooth boundary, and what we need is a control on the 
angles of the corners - where the boundaries of two or more such sets intersect.
In the non-uniform case we used the presence of many unipotents, and the nice
actions of unipotent groups on the boundary at infinity. In the compact case
there are no unipotents, so different tools are required.

Let $A,B\subset S$ be convex bodies with smooth boundary and with a common
interior point in their intersection. The angle between the boundaries at a 
common point
$x\in \partial A \cap \partial B$ is measured by $\pi$ minus the angle between
the external normal vectors $\hat n_A(x), \hat n_B(x)$
$$
 \phi_x(\partial A,\partial B)=
 \pi -\angle \big(\hat n_A(x),\hat n_B(x)\big).
$$ 
Thus, a big angle between the boundaries corresponds to a small angle between 
the normals.
The following lemma states that when $A$ and $B$ are sub-level sets 
of commuting isometries, these angles are $\geq\frac{\pi}{2}$. A similar statement appeared in \cite{BGS}.

For an isometry $\gc$, and for $x\in S$ with 
$d_{\gc}(x)=a>\min d_\gc$ 
we denote by $\hat n_{\gc}(x)$ the external (with respect to 
$\{ d_{\gc} \leq a \}$)
normal to 
$$
 \{ d_{\gc} =a \}=\partial\{ d_{\gc}\leq a \}.
$$

\begin{lem}\label{bigA}{\bf (Commutativity implies big angles):}
If the isometries $\alpha$ and
$\beta$ commute then $\hat n_{\alpha}(x)\cdot \hat n_{\beta}(x)\geq 0$.
\end{lem}

\begin{proof}
Let $a_\alpha =d_\alpha (x),a_\beta =d_\beta (x)$.
Since $\{ d_\alpha =a_\alpha \}$ is a level set for $d_\alpha$, we see that
$\hat n_{\alpha}(x)$ is the direction of the gradient
$(\nabla \cdot d_\alpha) (x)$.
Thus if $c(t)$ is the geodesic line through $x$ with 
$\dot c(0)=\hat n_{\beta}(x)$
then it is enough to show that
$$
 (\nabla \cdot d_\alpha )(x) \cdot \hat n_{\beta}(x)=
 \frac{d}{dt}|_{t=0}\{ d_\alpha \big( c(t)\big)\} \geq 0.
$$
Let $p$ denote the projection on the convex set 
$\{ d_{\beta} \leq a_\beta \}$. Since $\alpha$ and $\beta$ commute, the 
set $\{ d_{\beta} \leq a_\beta \}$ is $\alpha$-invariant and therefore
$\alpha$ commute with $p$. Since $p$ decreases distances and since the 
geodesic lines $c(t),\ga\cdot c(t)$ are both orthogonal to 
$\partial\{ d_\gb\leq a_\gb\}$ one sees that $d\big( c(t),\ga\cdot c(t)\big)$ 
is a non-decreasing function of $t$. Thus 
$$
 (\nabla \cdot d_\alpha )(x)\cdot \hat n_{\beta}(x)=
 \frac{d}{dt}|_{t=0}\{ d_\alpha \big( c(t)\big)\} =
 \frac{d}{dt}|_{t=0}\{d\big( c(t),\alpha \cdot c(t)\big)\} \geq 0.
$$
\end{proof}

We shall need similar information when $A,B$ are replaced by their 
$\gep$-neighborhoods $(A)_\gep ,(B)_\gep$. The following lemma explains that
the situation then only improves.
 
Define
$$
 \varphi_t= \sup_{x\in \partial (A)_t \cap \partial (B)_t}
 \angle\big(\hat n_{(A)_t}(x),\hat n_{(B)_t}(x)\big),
$$
then we have

\begin{lem}\label{monoton}{\bf (Monotonicity of angles):}
$\varphi_t$ is a non-increasing function of $t$.
\end{lem}

\begin{proof}
We need to show that if $t_1>t_2\geq 0$ then $\varphi_{t_1}\leq 
\varphi_{t_2}$. If $\varphi_t$ vanishes at some point $t=t_0$ then, as it is
easy to verify, the union $A_{t_0}\cup B_{t_0}$ is convex, and $\varphi_t=0$
for any $t>t_0$. We may therefore assume that $\varphi_t$ is strictly positive
in our segment $[t_2,t_1]$. 

Fix a common interior point 
$$
 y\in\textrm{int}(A)\cap\textrm{int}(B).
$$
Now 
$$
 \varphi_t=\sup_{x\in \partial (A)_t \cap \partial (B)_t}
 \angle \big(\hat n_{(A)_t}(x),\hat n_{(B)_t}(x)\big) =
 \sup_{R\subset S}\{\max_{x\in R\cap\partial (A)_t \cap \partial (B)_t} 
 \angle \big(\hat n_{(A)_t}(x),\hat n_{(B)_t}(x)\big)\}
$$
where $R$ runs over the (compact) balls centered at $y$.
It is therefore enough to show that 
$$
 \varphi_{t_2}\geq\max_{x\in R\cap\partial (A)_{t_1} \cap \partial (B)_{t_1}} 
 \angle \big(\hat n_{(A)_{t_1}}(x),\hat n_{(B)_{t_1}}(x)\big)
$$
for any such $R$.

Fix $R$ large enough (so that the intersection 
$R\cap\partial (A)_{t_1} \cap \partial (B)_{t_1}$ is not empty)
and let 
$$
 \tilde\varphi_t=\max_{x\in R\cap\partial (A)_t \cap 
 \partial (B)_t} \angle \big(\hat n_{(A)_t}(x),\hat n_{(B)_t}(x)\big).
$$
Since $\tilde\varphi_t$ is obviously continuous,
it is enough to show that for any $t_2 < t \leq t_1$ we have 
$\tilde\varphi_{t-\gD t}\geq\tilde\varphi_t$ whenever $\gD t$ is small enough.
Fix $t$, and let $\gD t$ be sufficiently small so that the argument below holds.

To simplify notation we replace $A$ (resp. $B$) by $(A)_{t-\gD t}$ 
(resp. $(B)_{t-\gD t}$) and assume that $t-\gD t=0$ (this is just a matter of
changing names after $t$ and $\gD t$ are fixed). 
Let $x_t\in R\cap\partial (A)_t \cap\partial (B)_t$ be such that 
$$
 \angle \big(\hat n_{(A)_t}(x_t),\hat n_{(B)_t}(x_t)\big) =\tilde\varphi_t.
$$
It is enough to show that there is $x_0\in R\cap\partial A \cap \partial B$
with
$$
 \angle \big(\hat n_A(x_0),\hat n_B(x_0)\big)\geq
 \angle \big(\hat n_{(A)_t}(x_t),\hat n_{(B)_t}(x_t)\big).
$$

Let $x_0\in \partial A \cap \partial B$ be a point at the minimal possible distance from $x_t$.
Denote by $\hat u_0$ (resp. $\hat u_t$) the tangent to the geodesic line
$\overline{x_0,x_t}$ at $x_0$ (resp. at $x_t$).

Since $\gD t$ is arbitrarily small, $\hat u_t$ is roughly in the direction of
the bisector of the angle between $\hat n_{(A)_t}(x_t)$ and 
$\hat n_{(B)_t}(x_t)$, and $x_0$ is closer to $y$ than $x_t$. In particular 
$x_0\in R$.

By the Lagrange multipliers theorem $\hat u_0$ is a linear combination of
$\hat n_A(x_0)$ and $\hat n_B(x_0)$ and since $\gD t$ is arbitrarily small, we
can assume that $\hat u_0$ is in the convex cone spanned by   
$\hat n_A(x_0)$ and $\hat n_B(x_0)$ (again, in the limit case $\hat u_0$
is the direction of the bisector of the angle between
$\hat n_A(x_0),\hat n_B(x_0)$). Thus 
$$ 
 \angle \big(\hat n_A(x_0),\hat n_B(x_0)\big) =
 \angle \big(\hat n_A(x_0),\hat u_0\big) +
 \angle \big(\hat u_0,\hat n_B(x_0)\big),
$$
and since (by the triangle inequality for angles) 
$$
 \angle \big(\hat n_{(A)_t}(x_t),\hat n_{(B)_t}(x_t)\big)\leq
 \angle \big(\hat n_{(A)_t}(x_t),\hat u_t\big)+
 \angle \big(\hat u_t,\hat n_{(B)_t}(x_t)\big),
$$
it is enough to show that 
$$
 \angle \big(\hat n_A(x_0),\hat u_0\big)\geq 
 \angle \big(\hat n_{(A)_t}(x_t),\hat u_t\big)
$$
(and the analogous inequality for $B$ instead of $A$ whose proof is the same).

Let $c(s)$ be the geodesic line of unit speed 
with $c(0)=x_0,\dot c(0)=\hat u_0$, and let
$c(s_0)=x_t$ (i.e. $s_0=d(x_0,x_t)$). Then the above inequality follows from
the following
$$
 \hat u_0\cdot \hat n_A(x_0)=\frac{d}{ds}|_{s=0}D_A\big( c(s)\big)\leq 
 \frac{d}{ds}|_{s=s_0}D_A\big( c(s)\big)=
 \hat u_t \cdot \hat n_{(A)_t}(x_t)
$$
and this follows from the convexity of the function $D_A$.
\end{proof}

The following lemma, explains how we can use the information obtained above,
in order to find a direction with respect to which the directional derivatives
of all the corresponding distance functions are large.

\begin{lem}{\bf (Existence of a good direction):}\label{b(d)}
There is a constant $b(d)>0$ such that for any set of unit vectors
$\{ \hat n_i\}_{i\in I}\subset \BR^d$, that satisfy the condition 
$\hat n_i\cdot \hat n_j \geq 0$ for any $i,j\in I$, there is a unit vector
$\hat f$ such that $\hat f\cdot \hat n_i > \frac{1}{b(d)}$ for any $i\in I$.
\end{lem}

\begin{proof}
Let $\gD$ be a maximal $1$-discrete subset of $\{ \hat n_i\}_{i\in I}$.
Then there are at most $\frac{b(d)}{2}$ elements in $\gD$, where $b(d)$ is 
some constant.
For any $\hat n_i \in \{ \hat n_i\}_{i\in I}$ there is 
$\hat n_{i_0}\in \gD$ with $\| \hat n_i-\hat n_{i_0}\| <1$, which implies
$\hat n_i \cdot \hat n_{i_0} >1/2$.
Let 
$$
 \hat f=\frac{\Sigma_{\hat n_j\in \gD} \hat n_j}{\|\Sigma_{\hat n_j\in \gD} 
 \hat n_j\|},
$$
then 
$$
 \hat n_i \cdot \hat f >\hat n_i \cdot \frac{\Sigma_{\hat n_j\in \gD} 
 \hat n_j}{\frac{b(d)}{2}}\geq
 \frac{2}{b(d)}\hat n_i \cdot \hat n_{i_0}\geq \frac{1}{b(d)}.
$$
\end{proof}

Summarizing the above discussion together with the arguments of the previous 
section, we conclude the following less general but more practical version of 
Lemma \ref{simplicial-complex}:

\begin{thm}\label{nice-boundary}
Let $S$ be a symmetric space of non-compact type, $M=\gC\backslash S$ an $S$-manifold, and
$\gep'\geq\gep >0$. Assume that: 
\begin{itemize}
\item 
$M'\subset M_{\geq\gep}$ is a submanifold with boundary
for which the complement of the pre-image in the universal covering is given
by 
$$
 \ti X=S\setminus\ti M'=\cup_{\gc\in\gC'}\{d_\gc <\gep'\},
$$
where $\gC'$ is a subset of $\gC$ which is invariant under conjugation by elements of $\gC$.
\item
$M'$ is homotopically equivalent to its $\frac{\gep}{2}$-shrinking
$)M'(_\frac{\gep}{2}$.
\item 
For any $x\in S$ the group 
$$
 \langle \gc\in\gC':d_\gc (x) \leq 3\gep'\rangle
$$
is either unipotent or commutative.
\end{itemize} 
Then there are positive constants $d,\ga$ which depend on $S,~\gep,\gep'$, 
such that $M'$ is homotopically equivalent to some 
$\big( d,\ga\cdot\vol (M')\big)$-simplicial complex.
\end{thm}

\begin{proof}
We shall show that the conditions of Lemma \ref{simplicial-complex} are 
satisfied. 
The sets $\{d_\gc <\gep'\}$ are convex and open with smooth boundary and 
$X$ is the locally finite union of them.\\

It is given that $M'$ is contained in the $\gep$-thick part and that it is
homotopic to its $\frac{\gep}{2}$-shrinking.\\

We should indicate how to define the constant $b$ and the direction 
$\hat n(x)$, from the conditions of Lemma \ref{simplicial-complex}.

We define
$$
 b=\max \{\frac{2}{\gb (\gep')},b(d)\}
$$ 
where $\gb$ is the function defined in the proof of Theorem \ref{thm1} in the 
previous section, and $b(d)$ is the
constant defined in Lemma \ref{b(d)} for $d=\dim S$.

Next we define the directions $\hat n(x)$. Let $x\in S\setminus\ti X$ be a 
point with $d(x,\ti X)\leq\gep$. If the group 
$\langle \gc\in\gC':d_\gc (x) \leq 3\gep'\rangle$ is unipotent, then we define
the direction $\hat n(x)$ in the same way as it is done in the previous 
section.
(It is shown in the previous section that in this case the conditions of 
Lemma \ref{simplicial-complex} are satisfied.)
Assume that the group $\langle \gc\in\gC':d_\gc (x) \leq 3\gep'\rangle$ is 
abelian. Let $A(x)\subset \gC'$ be the set of elements $\gc\in \gC'$
with $D_{\{d_\gc <\gep'\}}(x)=d(x,\ti X)$. Then for $\gc\in A(x)$
$$
 d_\gc (x)\leq 2d(x,\ti X)+\gep'\leq 3\gep',
$$
and hence $A(x)$ is contained in 
$\langle \gc\in\gC':d_\gc (x) \leq 3\gep'\rangle$. Thus $A(x)$ is abelian.
By Lemma \ref{bigA}, for any $\ga ,\gb\in A(x)$ and any 
$y\in\partial \{d_\ga <\gep'\}\cap\partial \{d_\gb <\gep'\}$, the inner 
product of the corresponding external normal vectors at $y$ is non-negative. 
By Lemma \ref{monoton},
the inner product at $T_x(S)$ of the external normal vectors to
$\partial (\{d_\ga <\gep'\})_{d(x,\ti X)}$ and
$\partial (\{d_\gb <\gep'\})_{d(x,\ti X)}$ at $x$, is also non-negative. 
Therefore it follows from Lemma \ref{b(d)} there is a direction $\hat f\in T_x(S)$ 
for which the inner product of $\hat f$ with the external normal vector at $x$ to
$\partial (\{d_\gc <\gep'\})_{d(x,\ti X)}$ is $\geq\frac{1}{b(d)}$ for any 
$\gc\in A(x)$.
Thus, we can take $\hat n (x)=\hat f$.
\end{proof}


\section{The proof of \ref{thmA}(2) in the rank-1 case}\label{rk-1}

Theorem \ref{thmA}(2) can be proved, independently of the rank, using the argument
that we shall present in section \ref{pfBC}. 
For two reasons we chose to do the rank-1
case separately. The first reason is that in the remaining higher rank case all
locally symmetric manifolds of finite volume are arithmetic. Since Theorem
\ref{thmA} already took care of the non-compact arithmetic case, we can assume compactness when 
considering higher rank manifolds. This would make the argument in section \ref{9} simpler.
The second reason is that in the rank-1 case there is another proof, which is
in some sense more simple. This proof, which we shall present below, is 
basically the same as the one given in \cite{BGLM} for the hyperbolic case, and with the tools
developed in sections 3,4,6 and 7, it could be applied to general rank-1 symmetric spaces which
are not necessarily of constant curvature.

\begin{rem}
For $S=\BH^2$ the statement of Conjecture \ref{conjA} follows 
from the Gauss Bonnet theorem without the arithmeticity assumption.
This is because the volume determines the genus and bound the possible number of
cusps. For this reason we allow 
ourselves to ignore this case in this and in the following section.
\end{rem}

Recall what we intend to prove:

\begin{thm}[{\bf Theorem \ref{thmA}(2) for the rank-1 case}]\label{rk-1thm}
If $S$ is a rank-$1$ symmetric space of dimension $\geq 4$, then for some 
constants $\ga ,d$, the fundamental group $\pi_1(M)$ of any 
$S$-manifold $M=\gC\backslash S$ is isomorphic to the fundamental group of some
$\big( d,\ga\cdot\vol (M)\big)$-simplicial complex.
\end{thm}

We shall use the ordinary thick-thin decomposition. 
Let $\gep_s$ be the constant from the Margulis lemma, and let 
$\gep =\frac{\gep_s}{2}$.

\begin{thm}[{\bf Thick-thin decomposition in rank-1}](see \cite{Th} section 4.5, and \cite{BGS} 
section 10)\label{thick-thin} 
Assume that $\rk (S)=1$. Let $M$ be an $S$-manifold of finite volume and $M_{\leq \epsilon}^0$ an 
arbitrary connected component of $M_{\leq \epsilon}$. Then $M_{\leq \epsilon}^0$ belongs to either 
one of the following two types:
\begin{enumerate}
\item 
A tube, i.e. a tubular neighborhood of a short geodesic. In this case 
$M_{\leq \epsilon}^0$ is topologically a ball-bundle 
over the circle. Its fundamental group $\pi_1(M_{\leq \epsilon}^0)$ is 
infinite cyclic, and in particular abelian.

or
 
\item  
A cusp. In which case $M_{\leq \epsilon}^0$ is homeomorphic to 
$\BR^{\geq 0}\times \partial M_{\leq \epsilon}^0$, where 
$\partial M_{\leq \epsilon}^0$ is topologically a sub-manifold of 
codimension $1$.
A connected component of its pre-image  $\tilde  M_{\leq \epsilon}^0$ is
given by $\tilde  M_{\leq \epsilon}^0=
\cup_{\gc \in \gC_0}\{ d_\gc \leq \epsilon \}$ where $\gC_0\leq\gC$ is a 
subgroup isomorphic to $\pi_1(M_{\leq\gep}^0)$, and it is a ``star-shaped''
neighborhood of some point $z\in S(\infty )$ (i.e. each geodesic line with
$c(\infty )=z$ enters once into $\tilde M_{\leq \epsilon}^0$ and stays in it).
$\gC_0$ consists of unipotent elements only,
and each $\gc \in \gC_0$
fixes $z$ and leaves the horospheres around $z$ invariant.
Moreover $\gC_0$ is metabelian.
\end{enumerate}
\end{thm}

There are only finitely many connected components of $M_{\leq\gep}$.
Moreover, as $\dim (M)\geq 3$ the boundary of each 
connected component of $M_{\leq\gep}$ is connected, and hence the thick 
part  $M_{\geq\gep}$ is connected. Each cusp is homotopically equivalent to its
boundary. As $\dim (M)\geq 4$ each tube is a ball bundle over the circle,
with fibers of dimension $\geq 3$, and hence, its boundary is a sphere bundle
for a sphere of dimension $\geq 2$. Therefore, for each connected component of
the thin part, the injection of the boundary into the component induces an
isomorphism between the fundamental groups (which are both $\cong\BZ$). 
Thus by Van-Kampen's theorem we obtain:

\begin{cor}
$$
 \pi_1(M)\cong\pi_1(M_{\geq\gep}).
$$
\end{cor}

Take $M'=M_{\geq\gep}$. Then the following claim 
finishes the proof of Theorem \ref{rk-1thm}.

\begin{clm}
$M'$ satisfies the condition of Theorem \ref{nice-boundary}, with respect 
to $\gep'=\gep$.
\end{clm}

\begin{proof}
We only need to explain why $M'$ is homotopic to its 
$\frac{\gep}{2}$-shrinking, since the other conditions of Theorem
\ref{nice-boundary} follow directly from the definition of $M'$ and from Theorem \ref{thick-thin}. 

Since $\gep =\frac{\gep_s}{2}$, the $\frac{\gep}{2}$-neighborhoods of different 
connected components of
$M_{\leq\gep}$ are still disjoint, and hence, we can prove the homotopy 
equivalence by
showing that there is a deformation retract from 
$\overline{(M_{\leq\gep}^0)_\frac{\gep}{2}}\setminus M_{\leq\gep}^0$ to
$\partial (M_{\leq\gep}^0)_\frac{\gep}{2}$, for each connected component 
$M_{\leq\gep}^0$ of $M_{\leq\gep}$.

If $M_{\leq\gep}^0$ is a cusp, then one can define such a deformation retract 
by letting each point flow (at constant speed = the initial distance) along the 
unique geodesic line which connects it to the unique end of the cusp, i.e. along
the geodesic whose lifting in the universal covering converges, when $t\to -\infty$, to the unique
limit point in $S(\infty )$ of a lifting of the cusp.

If $M_{\leq\gep}^0$ is a tube, then the lifting of $M_{\leq\gep}^0$ and 
its $\frac{\gep}{2}$-neighborhood in $S=\ti M$ are both star-shaped with 
respect to the lifting $c$ of the short closed geodesic which lies inside $M_{\leq\gep}^0$. 
The generalized star-contraction (see \ref{GSC}) from 
$\overline{(\ti M_{\leq\gep}^0)_\frac{\gep}{2}}\setminus \ti M_{\leq\gep}^0$ to
$\partial (\ti M_{\leq\gep}^0)_\frac{\gep}{2}$ 
projects to a deformation retract of the corresponding subsets of $M$. 
\end{proof}


\section{The proofs of \ref{thmA}(2) and \ref{thmA}(3)}\label{pfBC}

We shall use the fact that when a low dimensional submanifold is 
removed from a high dimensional manifold, it doesn't change the low dimensional
homotopy groups. More precisely:

\begin{lem}\label{M-N}
Let $M$ be a connected manifold and $N\subset M$ a closed (not 
necessarily connected) sub-manifold.
\begin{itemize}
\item If $\text{codim}_M(N)\geq 2$ then there is a surjective homomorphism
$$
 \pi_1(M\setminus N)\to\pi_1(M).
$$
\item If $\text{codim}_M(N)\geq 3$ then $\pi_1(M\setminus N)\cong\pi_1(M)$.
\end{itemize}
\end{lem}

\begin{proof}
If $\text{codim}_M(N)\geq 2$ then any closed loop can be pushed to one which 
dose not intersect $N$. Similarly, if $\text{codim}_M(N)\geq 3$ then any 
homotopy of loops can be pushed to $M\setminus N$.
\end{proof}

The submanifold we are about to remove consists
of the union of all short closed geodesics of some certain type. Observe that
the lifting to $S=\ti M$ of a closed geodesic in $M=\gC\backslash S$ which corresponds to
$\gc\in\gC$ is an axis of $\gc$. In particular, any two closed geodesics from
the same homotopy class are parallel.

\begin{lem}\label{codimN}
Let $\gc\in G^0$ be a hyperbolic isometry which projects non-trivially to each 
simple factor of $G^0$, and let $\ti N\subset S$ be the union of all geodesics 
which are axes of $\gc$.

\begin{itemize}
\item If $S$ is not isometric to $\BH^2$ then $\text{codim}_S(\ti N)\geq 2$.
\item If additionally $S$ is neither isometric to $\BH^3, \BH^2\times\BH^2$ nor to
$\PSL_3(\BR )/\PSO_3(\BR )$, then $\text{codim}_S(\ti N)\geq 3$.
\end{itemize}
\end{lem}

\begin{proof}
Let $S^*$ be an irreducible factor of $S$. Denote by $\gc^*$ the projection of 
$\gc$ to the corresponding simple factor $G^*$ of $G^0$.
We need to estimate $\text{codim}_{S^*}\big(\min (\gc^*)\big)$. 

If $\gc^*$ fixes a point then we just remark that since the projection $\gc^*$
is non-trivial, 
$\text{codim}_{S^*}\big(\min (\gc^*)\big) >0$. Otherwise let $c$ be an axis
of $\gc^*$ in $S^*$ and let
$$
 G^*(c) = \{ g\in G^*:g\cdot c \text{~is parallel to~} c\}.
$$
An alternative way to define $G^*(c)$ is as follows. Pick two points
$p,q\in c(\BR )$ and let $g_{p,q}=\gs_p\cdot\gs_q$ be the corresponding 
transvection,  then $G^*(c)$ coincides with the centralizer group 
$C_{G^*}(g_{p,q})\leq G^*$ of $g_{p,q}$.
$G^*(c)$ is a closed reductive Lie subgroup of $G^*$ which acts 
transitively on $N^*=$ the union of all geodesics parallel to $c$ in $S^*$. 
Observe that $N^*$ contains the projection of $\ti N$ to $S^*$.

Choose an Iwasawa decomposition of $G^*$ which induces an Iwasawa 
decomposition of $G^*(c)^0$, i.e. write the Iwasawa decompositions of
$G^*$ and $G^*(c)$ simultaneously 
$$
 G^*=K\cdot A\cdot U, \text{~and~} G^*(c)^0=K(c)\cdot A\cdot U(c),
$$
such that $K(c)=K\cap G^*(c)^0$ and $U(c)=U\cap G^*(c)^0$.
(This could be done by choosing a flat $F\supset c$ and defining the torus $A$
to be the subgroup of $G^*$ which acts on $F$ by translations, and then 
choosing an order on the dual of $\mathfrak{a}=\text{Lie}(A)$ such that if 
$c(t)$ is given by $c(t)=e^{tH}$ for $H\in\mathfrak{a}$, then $\ga (H)\geq 0$ 
for any simple 
root $\ga$, and then taking the Iwasawa decomposition which corresponds to 
$A$ with this ordering.) 

Then $\dim (S^*)=\dim (A\cdot U)$ while $\dim (N^*)=\dim \big( A\cdot U(c)\big)$,
and hence
$$
 \text{codim}_{S^*}(N^*)=\text{codim}_U\big( U(c)\big).
$$
Let $P_{c(-\infty )}$ be the parabolic subgroup of $G^*$ which corresponds to
the point $c(-\infty )\in S^*(\infty )$, and let $P^-$ be the minimal parabolic
opposite to $P^+=AU$. Then, as it is easy to verify, the Lie algebra of $G^*$
is the direct sum
$$
 \text{Lie}(G^*)=\text{Lie}(U)+\text{Lie}(P^-),
$$
while the Lie algebra of the parabolic $P_{c(-\infty )}$ is the direct sum
$$
 \text{Lie}\big( P_{c(-\infty )}\big)=
 \text{Lie}\big( U(c)\big)+\text{Lie}(P^-).
$$
Hence $\text{codim}_U\big( U(c)\big)=\text{codim}_{G^*}(P_{c(-\infty )})$.
The lemma follows from the following proposition by a standard
case analysis.
  
\begin{prop}[see \cite{BaNe} lemma 3.4 and corollary 8]
Let $H$ be a connected simple Lie group and $P\leq H$ a proper closed connected 
subgroup. Then $\text{codim}_H(P)\geq\rk (H)$, and the equality can hold only
if $H$ is locally isomorphic to $\SL_n(\BR )$.
\end{prop}
\end{proof}

We shall now prove: 

\begin{thm}[{\bf Theorem \ref{thmA}(2) from the introduction
in the general case}]\label{thm2}
Let $S$ be a symmetric space of non-compact type, not isometric to
$\BH^2,\BH^3, \BH^2\times\BH^2, \PSL_3(\BR )/\PSO_3(\BR )$.
Then there are constants $\ga ,d$ such that the fundamental group of any
irreducible $S$-manifold of finite volume $M$ is isomorphic to the fundamental 
group of some $\big( d,\ga\cdot\vol (M)\big)$-simplicial complex.
\end{thm}

\begin{proof}[Proof of Theorem \ref{thm2}] 
Let $G=\text{Isom}(S)$, let $r$ be the rank of $S$, i.e. the real rank of $G$, let $n_s$ be the
index from the Margulis lemma and let $n=n_s!$. Fix $\gep$ to be one third of 
the $\gep_s$ from the Margulis lemma.

Since we already proved the theorem in the rank one case, we can assume that $\rk (S)\geq 2$.
Then any $S$-manifold is arithmetic. Moreover, as Theorem \ref{thm1} implies Theorem \ref{thm2}
for non-compact arithmetic manifold, we may assume that $M=\gC\backslash S$ is compact.
In this case, we have the following strengthening of the Margulis lemma.

\begin{lem}
For any $x\in S$, the group $\gC_{\gep_s} (x)$ contains an abelian subgroup of
index $n_s$.
\end{lem}

\begin{proof}
By the Margulis lemma $\gC_{\gep_s} (x)$ contains a subgroup 
$\gC_{\gep_s}^0 (x)$ of index $n_s$ which is contained in a connected nilpotent
Lie subgroup of $G$. By Lie's theorem $\gC_{\gep_s}^0 (x)$ is triangulable
over $\BC$. Therefore its commutator $[\gC_{\gep_s}^0 (x),\gC_{\gep_s}^0 (x)]$
contains only unipotents. As $\gC$ is cocompact it has no 
non-trivial unipotents. Thus $\gC_{\gep_s}^0 (x)$ is abelian.
\end{proof}

It follows that when $\gc_1 ,\gc_2\in\gC$ satisfy 
$\{ d_{\gc_1}\leq\gep\}\cap\{ d_{\gc_2}\leq\gep\}\neq\emptyset$ then $\gc_1^n$
and $\gc_2^n$ commute.

Consider an element $\gc\in\gC$. Replacing $\gc$ by $\gc^{[G:G^0]}$, we may
assume that $\gc\in G^0$. As $\gC\cap G^0$ is a uniform 
lattice in $G^0$, $\gc$ is semisimple.
Consider the sequence of centralizers
$$
 C_G(\gc)\leq C_G(\gc^n)\leq C_G(\gc^{n^2})\leq C_G(\gc^{n^3})\leq\ldots
 \leq C_G(\gc^{n^r}).
$$
As the centralizer of a semisimple element is determined by the type of the 
singularity of the element and by a torus which contains it, two consecutive
terms $C_G(\gc^{n^i})$ and $C_G(\gc^{n^{i+1}})$ in this sequence 
must coincide. Take $i$ to be the first time at which this happens and write 
$\gc'=\gc^{n^i}$. In this way we attach $\gc'$ to each $\gc\in\gC$.
Notice that if $\gc\in G^0$ then $\gc'=\gc^j$ for some $j\leq n^r$, 
and in general $\gc'=\gc^j$ for some $j\leq [G:G^0]n^r$. Set $m=[G:G^0]n^r$.

The reason we prefer to work with the $\gc'$'s is the following. If 
$$
\{d_{\gc_1'}\leq 3\gep\}\cap\{d_{\gc_2'}\leq 3\gep\}\neq\emptyset
$$ 
then $\gc_1'$ and $\gc_2'$ commute. Of course, by the above lemma, the 
non-empty 
intersection implies that $\gc_1'^n,\gc_2'^n$ commute. Since $\gc_1'$ has the 
same centralizer as $\gc_1'^n$, this implies that $\gc_1',\gc_2'^n$ commute, 
and since $\gc_2'$ has the same centralizer as $\gc_2'^n$,  
this implies that $\gc_1',\gc_2'$ commute.

Let $N\subset M$ be the subset which consists of the union of all
closed geodesics of length $\leq\gep$ which correspond to elements of the
form $\gc'$ for $\gc\neq 1$. In other words, its pre-image in $S$ is given by 
$$
 \ti N=\cup\{\min (\gc' ):\gc\in\gC\setminus\{ 1\} ,\min (d_{\gc'})\leq\gep\}.
$$ 
Then $N$ is a finite union of totally geodesic closed submanifolds. 
Since $\gC$ is irreducible and $G$ is center free, any $\gc\in\gC\setminus\{ 1\}$ 
projects non-trivially to each factor of $G$. So it follows from
Lemma \ref{codimN} that $\text{codim}_M(N)\geq 3$, and hence, by Lemma 
\ref{M-N}
$$
 \gC\cong\pi_1(M)\cong\pi_1(M\setminus N).
$$

Define
$$
 \ti X=\cup_{\gc\in\gC\setminus\{ 1\}}\{ d_{\gc'}<\gep\},~~
 \ti M'=S\setminus \ti X,
$$
and 
$$
 X=\gC\backslash\ti X,~~M'=\gC\backslash\ti M'=M\setminus X.
$$

Clearly $N\subset X$.
In order to prove the theorem we shall show:
\begin{enumerate}
\item 
$M\setminus N$ and $M'$ are homotopically equivalent.
\item 
$M'$ satisfies the conditions of Theorem \ref{nice-boundary}
with respect to $\gep$ and $\frac{\gep}{m}$ (corresponding to $\gep'$ and 
$\gep$ respectively in \ref{nice-boundary} ).
\end{enumerate}

{\bf Proof of (1):}
We shall construct the desired homotopy in two steps.
Since $M$ and $N$ are compact, and each connected component of $N$ is a 
finite union of totally geodesic submanifolds, there is small positive number 
$\eta>0$, such that the 
$\eta$-neighborhoods of the components of $N$ are still disjoint and contained 
in $X$. It is easy to verify that if $\eta$ is small,
$M\setminus N$ and $M\setminus \overline{(N)_\eta}$ are diffeomorphic.

Next, we claim that there is a deformation retract from 
$M\setminus \overline{(N)_\eta}$ to $M'$. In order to show this, we shall
apply Lemma \ref{deformation-retract} to the 
set of functions 
$$
 \mathcal{F}=\{d_{\gc'}:\gc\in\gC\setminus\{1\},\min (d_{\gc'})<\gep\}
$$ 
on $Y=L=S\setminus\overline{(\ti N)_\eta}$, and $M'=\mathcal{F}_{\geq\gep}$. 
We have to indicate what are the directions 
$\hat n(x)$, and what is the continuous function $\gb$.

Let $x\in \ti X\setminus\overline{(\ti N)_\eta}$, and let 
$\{ d_{\gc'_1},\ldots ,d_{\gc'_k}\}\subset \mathcal{F}$ be the subset of
functions satisfying $d_{\gc'_i}(x)\leq 3\gep$. Then the group
$\langle \gc'_1,\ldots ,\gc'_k\rangle$ is abelian. Thus for each $j\leq k$
the convex set $\cap_{i=1}^{j-1}\min (\gc'_i)$ is $\gc'_j$ invariant, 
and hence, by induction 
$$
\cap_{i=1}^j\min (\gc'_i)\neq\emptyset.
$$
Pick arbitrarily $y\in\cap_{i=1}^k\min ({\gc'_i})$ and define 
$\hat n(x)$ to be the tangent at $x$ 
to the geodesic line which goes from $y$ to $x$.

For $t\leq 3\gep$ we define $\gb (t)$ to be the minimum of the directional 
derivative of $d_{\gc'}$ at $x$ with respect to the tangent to the geodesics 
$\overline{z,x}$, where this minimum is taken over all 
$x\in\ti X\setminus (\ti N)_\eta$, all $d_{\gc'}\in \mathcal{F}$ with 
$d_{\gc'}(x)=t$, and all
$z\in\min (\gc' )$.
Since 
\begin{itemize}
\item
$N$ and $M\setminus (N)_\eta$ are compact, 
\item
up to conjugations there are only finitely many $\gc'$'s with 
$\min (d_{\gc'})<\gep'$, and 
\item 
for any selection
of $x,z,\gc'$ as above, the corresponding directional derivative is positive,
\end{itemize}
it follows that $\gb$ is a well defined
continuous positive function.\\

{\bf Proof of (2):}
We shall check that the conditions of Theorem
\ref{nice-boundary} are satisfied.

Let $x\in S$ and assume $d_\gc (x)\leq\frac{\gep}{m}$
for some $\gc\neq 1$ in $\gC$. As $\gc'=\gc^j$ for some $j\leq m$, it follows
that 
$$
d_{\gc'}(x)=d_{\gc^j}(x)\leq j\cdot d_\gc (x)\leq\frac{j\cdot\gep}{m}\leq\gep.
$$
Thus $x\in\ti X$. This shows that $M'$ is contained in 
$M_{\geq\frac{\gep}{m}}$.

If $\{ d_{\gc_1'}\leq 3\gep\}\cap\{ d_{\gc_2'}\leq 3\gep\}$ then 
$\gc'_1$ commutes with $\gc'_2$. Thus, the last condition of Theorem
\ref{nice-boundary} is also satisfied.

We shall show that there is a deformation retract from $M'$ to its 
$\frac{\gep}{2m}$-shrinking, by an analogous way to the second step in the 
proof of (1) above. 
This time take
$$
 \mathcal{F}=\{D_{\{ d_{\gc'}\leq\gep\}}:\gc\in\gC\setminus\{1\},\min (d_{\gc'})<\gep\}
$$ 
while $\frac{\gep}{2m}$ plays the role of $\gep$ in Lemma 
\ref{deformation-retract}.
We define the directions $\hat n(x)\in T_x(S)$ 
to be the tangent to the geodesic $\overline{y,x}$ for arbitrary
$y\in\cap\{\min (\gc' ):d_{\gc'}\in\Psi_{x,3\gep}\}$. 
Additionally, we define
$\gb (t)$ to be the minimum of the directional 
derivative of $D_{\{d_{\gc'}\leq\gep\}}$ at $x$ with respect to the tangent 
of the geodesics $\overline{z,x}$, where the minimum is taken over all 
$x\in S$, all $D_{\{d_{\gc'}\leq\gep\}}\in \mathcal{F}$ with 
$D_{\{d_{\gc'}\leq\gep\}}(x)=t$, and all
$z\in\min (\gc' )$. Since for $t\geq 0$ and for a relevant $\gc'$, the set
$\{ D_{\{d_{\gc'}\leq\gep\}}\leq t\}$ is a convex body with smooth boundary and since
any $z$ as above belongs to the interior of this set, the directional
derivatives mentioned above are always positive. By compactness we get that
$\gb$ is a continuous positive function defined for any $0\leq t\leq\frac{3\gep}{2m}$.
This finishes the proof of Theorem \ref{thm2} 
\end{proof}

Next we prove:

\begin{thm}[{\bf Theorem \ref{thmA}(3) of the introduction}]\label{thm3}
For any $S$, there are constants $\ga ,d$, such that the fundamental group
of any irreducible $S$-manifold $M$ is isomorphic to a quotient of the fundamental group 
of some $\big( d,\ga\cdot\vol (M)\big)$-simplicial complex. 
\end{thm}

\begin{proof}
There are only $3$ cases left to deal with. $\PSL_3(\BR )/\PSO_3(\BR ),
\BH^2\times\BH^2$ and $\BH^3$.

The proof for the first two cases goes verbatim as the proof of Theorem 
\ref{thm2}, since these cases are of higher rank and we can assume compactness.
The only difference is that in these cases, we have only $\text{codim}_M(N)\geq 2$ (instead of $\geq 3$) by 
the first part of Lemma \ref{codimN}, 
so the result follows from the first part of Lemma \ref{M-N}.

For hyperbolic $3$-manifolds, one should only throw out finitely many circles
which are all the closed geodesics of length $\leq\gep$ and then prove that
what is left is homotopically equivalent to the $\gep$-thick part. This is 
done by deforming each cusp to its boundary as in section 
\ref{rk-1}, and deforming each compact component minus a circle to its 
boundary, as it is done in the proof of \ref{thm2}. 
Then, again, Theorem \ref{nice-boundary} finishes the proof. 
\end{proof}


\section{Some remarks on Conjecture \ref{conjA}, Theorem \ref{thmA} and their relations to algebraic number 
theory}\label{9} 

Let $p(x)\in\BZ [x]$ be an integral monic polynomial, and let
$$
 p(x)=\prod_{i=1}^k(x-\ga_i)
$$ 
be its factorization into linear factors over $\BC$. Denote by $m(p)$ its
{\bf exponential Mahler measure}
$$
 m(p)=\prod_{|\ga_i|>1}|\ga_i|.
$$
The following is known as Lehmer's conjecture. 

\begin{conj}
There exists a constant $\ell >0$ such that if $p(x)$ is an integral 
monic polynomial with $m(P)\neq 1$, then $m(p)>1+\ell$.
\end{conj}   

Denote by $d(p)$ the number of roots $\ga_i$ with absolute value $>1$
$$
 d(p)=\#\{ \ga_i:|\ga_i|>1\}.
$$

The following conjecture of Margulis 
is weaker than Lehmer's conjecture.

\begin{conj}[see \cite{Mar1} $IX$ 4.21]\label{MargulisConj}
There is a function $\ell :\BN\to\BR^{>0}$ such that 
$m(p)\geq 1+\ell \big( d(p)\big)$ for any non-cyclotomic monic
polynomial $p(x)\in\BZ [x]$.
\end{conj}

If Conjecture \ref{MargulisConj} is true, then for any 
symmetric space of  non-compact type $S$, the minimal injectivity radius of any compact 
arithmetic $S$-manifold is bounded from below by some positive constant $r=r(S)$ (see also 
\cite{Mar1}, page 322 for a similar statement). 
To see this, we argue as follows: 
Let $G^0$ be the identity component of $G=\textrm{Isom}(S)$. As $G^0$
is center-free we can identify it with its adjoint group 
$\textrm{Ad}(G^0)\leq\textrm{GL}(\mathfrak{g})$.
Let $\gC\leq G^0$ be a torsion-free uniform arithmetic lattice in $G^0$. 
We think of $\gC$ as the intersection of the fundamental group of some compact
arithmetic $S$-manifold with $G^0$.
Since $\gC$ is arithmetic, there is a compact extension 
$G^0\times O$ of $G^0$ and a $\BQ$-rational structure on the Lie algebra 
$\mathfrak{g}\times \mathfrak{o}$ of $G^0\times O$,
such that $\gC$ is the projection to $G^0$ of a lattice $\ti{\gC}$, which is 
contained in $(G^0\times O)_{\BQ}$ and commensurable to the group of integral
points $(G^0\times O)_{\BZ}$ with respect to some $\BQ$-base of
$(\mathfrak{g}\times \mathfrak{o})_{\BQ}$. By changing this $\BQ$-base, we 
can assume that $\ti{\gC}$ is in fact contained in $(G^0\times O)_{\BZ}$.
This means that the characteristic polynomial $p_{\ti{\gc}}$ of any 
$\ti{\gc}\in\ti{\gC}$ is a monic integral polynomial. 
As $\gC$ is discrete and
torsion-free, $m(p_{\ti{\gc}})>1$ for any $\ti{\gc}\in\ti{\gC}$
which projects to a non-trivial element in $\gC$.
Since $O$ is compact, any eigenvalue of $\tilde\gc$ with absolute value 
different from $1$ is also an eigenvalue of its projection $\gc\in G^0$. 
In particular
$$
 m(p_{\gc})= m(p_{\ti\gc})\geq 1+\min_{i\leq\dim G}\ell (i).
$$ 
We conclude that $\gc$ is outside the open set 
$$
 U=\{ g\in G^0: m(p_g)<1+\min_{i\leq\dim G}\ell (i)\}.
$$ 
Clearly $U$ contains any
compact subgroup of $G$, and contains a subset of the form
$$
 \{ g\in G^0:g \textrm{~is semisimple and~} \min d_g < \ti r(S)\}
$$ 
for some positive constant $\ti r(S)$.
Hence the minimal injectivity radius of $\gC\backslash S$ is $\geq \ti r(S)$.
Finally, if $\gC\leq G$ is a torsion free lattice which is not necessarily 
contained in $G^0$ then the minimal injectivity radius of $\gC\backslash S$ is at 
least $\frac{\ti r(S)}{[G:G^0]}=r(S)$.
This implies: 

\begin{cor}
Conjecture \ref{MargulisConj} implies Conjecture \ref{conjA} (for compact 
arithmetic manifolds).
\end{cor}

\begin{proof}
For a given $M$, choose a maximal $r$-discrete net $\mathcal{C}$, and take $\mathcal{R}$ to be the
simplicial complex which corresponds to the nerve of the cover of $M$ by the $r$-balls whose centers form $\mathcal{C}$.
Then $M$ is homotopic to $\mathcal{R}$ which is a $\big(\frac{\vol (B_{2.5r})}{\vol (B_{r/2})},
\frac{\vol (M)}{\vol (B_r)}\big)$-simplicial complex.
\end{proof}

For compact locally symmetric manifolds, the minimal injectivity radius  
equals half of the length of the shortest close geodesic. 
For non compact arithmetic $S$-manifolds we already know the absence of short 
closed geodesics (see Remark \ref{5.7}). Therefore
Conjecture \ref{MargulisConj} implies 

\begin{conj}\label{injrarCong}
For any $S$, there exists a constant $l=l(S)$ such that no arithmetic 
$S$-manifold contains a closed geodesic of length $\leq l$.
\end{conj}

A simple argument, which uses the fact that for non-cyclotomic monic integral
polynomials $F(x)$ of a fixed degree $k=\dim (G)+\dim (O)$, $m(F)$ is bounded
away from $1$, shows that 
for all the compact arithmetic $S$-manifolds, which
arise by constructions in which the compact extending group $O$ is fixed, 
the infimum on the minimal injectivity radius is positive, and the statement
of Conjecture \ref{injrarCong} holds, independently of the rational structure 
and of the specific choice of a manifold within a commensurability class. 
More precisely:

\begin{prop}\label{fixedO}
Given a compact semi-simple Lie group $O$, there is a positive constant
$r=r(S,O)$, such that the minimal injectivity radius is $\geq r$,
for any compact manifold $M=\gC\backslash S$ such that $\gC\cap G^0$ is commensurable to the group 
$\pi_{G^0}\big( (G^0\times O)_\BZ\big)$ for some $\BZ$-structure on 
$\mathfrak{g}\times \mathfrak{o}$.
In particular, any such manifold is homotopically equivalent to a 
$\big(\frac{\vol (B_{2.5r})}{\vol (B_{r/2})},
\frac{\vol (M)}{\vol (B_r)}\big)$-simplicial complex.
\end{prop}
 
A real algebraic integer $\gt >1$ is called a {\bf Salem number} if all its 
conjugates in $\BC$ have absolute value $\leq 1$.
In Conjecture \ref{MargulisConj} it is not even known whether $\ell (1)>0$, i.e.
whether there is a positive gap between $1$ and the set of Salem 
numbers. Sury \cite{sury} showed that the existence of such a gap is equivalent to the 
existence of an identity neighborhood in $\SL_2(\BR )$ which intersects 
trivially any uniform arithmetic lattice. 
Therefore the existence of such a gap implies
a positive infimum on the length of closed geodesics, when considering all the 
arithmetic surfaces.
We shall now discuss the relation between this gap and 
manifolds locally isometric to $S=\textrm{SL}_3(\BR )/\textrm{SO}_3(\BR )$ or to 
$\BH^2\times\BH^2$ (the two cases for which we could not prove the analog of Theorem \ref{thmA}(2)).

\begin{clm}
Assume there is a positive gap between $1$ and the set of Salem numbers (i.e. $\ell (1)>0$).
Then the analog of Theorem \ref{thmA}(2) holds also for the symmetric space
$S=\textrm{SL}_3(\BR )/\textrm{SO}_3(\BR )$. 
\end{clm}

\begin{proof}
Let $\gt >1$, and assume that $\gc =\textrm{diag}(\gt ,\gt ,\gt^{-2})$
is an element of a lattice $\gC\leq\textrm{SL}_3(\BR )$. One 
can easily compute the eigenvalues of $\textrm{Ad}(\gc )$ : they are 
$\gt^3,\gt^{-3}$ and $1$. As $\gC$ is arithmetic, $\gt$ is a Salem number, 
and thus, by
our assumption, is bounded away from $1$. We conclude that the displacement 
functions $d_\gc$ of such elements are uniformly bounded away from $0$.
Take 
$$
 \gep^* < \text{inf}\{\min d_\gc :\gc =\textrm{diag}(\gt ,\gt ,\gt^{-2}), \gt 
 \text{~is a Salem number}\}
$$
which is also smaller than $\gep_s$ from the Margulis lemma.

If $g$ is a hyperbolic element of an arithmetic lattice in $\SL_3(\BR )$ with $\min (d_g )\leq\gep^*$, then
$g$ has no real eigenvalues of multiplicity $2$. It is easy to see that $\min (g)$ is then a flat or a single
geodesic. In particular $\dim\big(\min (g)\big)\leq 2$.
As $\dim (S)=5$, the submanifold $N$ which consists of the union of all closed
geodesics of length $\leq\gep$ has codimension $\geq 3$. 
Therefore we can apply word by word the argument of
the proof of Theorem \ref{thm2} for $S$ with $\gep^*$.
\end{proof}     

If moreover $\ell (2)>0$, then the same is true also for the 
second remaining case $S=\BH^2\times\BH^2$. 

\begin{clm}
Assume $\ell (1)\cdot\ell (2)>0$. Then the analog of \ref{thmA}(2) holds also for $S=\BH^2\times\BH^2$.
\end{clm}

\begin{proof}
In contrast with the situation for $\SL_3(\BR )$, our problem here arises only for {\it regular}
elements $\gc\in\SL_2(\BR )\times\SL_2(\BR )$ (for otherwise
$\min (\gc )$ is a single geodesic and its codimension is $3$). 
But such an element has always
the form $\gc =(\gc_1,\gc_2)$ where $\gc_i\in\SL_2(\BR )$ are diagonalizable
over $\BR$. If $\gc_i$ is conjugate to 
$\text{diag}(\gt_i,\gt_i^{-1}),~\gt_i >1$
then the only possible non-trivial conjugate of $\gt_1$ outside the unit disk
is $\gt_2$.
Therefore, assuming $\ell (1),\ell (2)>0$ we can choose $\gep^*$ small enough, so that all 
arithmetic elements with minimal displacement $\leq\gep^*$ are singular, and each has a unique axis.
Then we can apply the argument of the proof of Theorem \ref{thm2} to this case. 
\end{proof}

\begin{rem}
Similarly, $\ell (1)>0$ implies also the analog of  
\ref{thmA}(2) for compact arithmetic $3$-manifolds.
\end{rem}


\section{Estimating the size of a minimal presentation}

If Conjecture \ref{conjA} is true, 
then, given a symmetric space of non-compact type $S$, and an $S$-manifold $M$, 
such that either
\begin{itemize}
\item $S$ is not $\BH^3$, or
\item $M$ is arithmetic,
\end{itemize}
the minimal size of a presentation for the fundamental group should be bounded
linearly by the volume. 
From Theorem \ref{thmA} we deduce this for most cases: 

\begin{defn}
We say that a presentation of a group $\langle \gS : W\rangle$ is 
{\bf standard} if the length of each $w\in W$ is $\leq 3$.
\end{defn}

\begin{thm}\label{presentation-thm}
Assume that either
\begin{itemize}
 \item $S$ is not isomorphic to $\BH^3, \BH^2\times\BH^2, \PSL_3(\BR )/\PSO_3(\BR )$, or
 \item $M$ is non-compact arithmetic.
\end{itemize}
Then for some constant $\eta =\eta (S)$, independent of $M$,
the fundamental group $\pi_1(M)$ admits a standard presentation 
$$
 \pi_1(M)\cong\langle \gS : W\rangle
$$
with $|\gS |,|W|\leq\eta\cdot\vol (M)$.

\end{thm}

\begin{proof}
Let $\mathcal{R}$ be the $\big( d,\ga\cdot\vol (M)\big)$-simplicial complex which corresponds to $M$ by Theorem
\ref{thmA} (1) or (2). Fix a spanning tree $T$ for 
$\mathcal{R}$, and take the generating set $\gS$ for $\pi_1(\mathcal{R})\cong\pi_1(M)$ 
which consists those closed loops which contain exactly one edge outside $T$.
We thus obtain a generating set of size less then the number of edges of the 
1-skeleton $\mathcal{R}^1$ which is at most $\frac{\ga (S)\vol (M)d(S)}{2}$.
In other words, we take for each edge of $\mathcal{R}^1\setminus T$ the element
of $\pi_1(\mathcal{R}^1)$ which corresponds to the unique cycle (with 
arbitrarily chosen orientation) which is obtained by adding this edge to $T$.
Additionally, let the set of relations $W$ consist exactly 
those words which are induced from 2-simplexes of $\mathcal{R}^2$ (we take
one such relation for each 2-simplex). In this way we obtain
a set of relations of size $\leq\ga (S)\vol (M) d^2(S)$ which is a bound for
the number of triangles in $\mathcal{R}^1$.
Thus, 
$$
 \eta =\ga (S)d^2(S)=\max\{ \ga (S)d^2(S),~\ga (S)d(S)/2\}
$$
will do.
Finally, the length of each $w\in W$ is exactly the number of edges in the 
corresponding 2-simplex which lie outside $T$, and thus the 
presentation $\langle \gS : W\rangle$ is standard.
\end{proof} 

\begin{rem}
Lower bounds for the size of any presentation, are known for hyperbolic 
$3$-manifolds (see \cite{cooper}). In this case, the upper bound obtained 
above (for non-compact arithmetic $3$-manifolds) is tight. 
\end{rem}

For non-arithmetic hyperbolic $3$-manifolds the analogous statement is 
evidently false (see Remark \ref{55}). Surprisingly the following result holds:

\begin{thm}\label{Pre-Hyp}
There is a constant $\eta$ such that the fundamental group of 
any complete hyperbolic $3$-manifold $M$ admits a presentation
$$
 \pi_1(M)\cong\langle \gS : W\rangle
$$
for which both $|\gS |$ and $|W|$ are $\leq\eta\cdot\vol (M)$.
\end{thm}

\begin{lem}
Let $S$ be a rank one symmetric space, and let $\gep =\frac{\gep_s}{3}$
where $\gep_s$ is the constant from the Margulis lemma. 
There is a constant $c=c(S)$ such that for every $S$-manifold $M$, the number of closed geodesics of length $\leq\gep$ in $M$ is at most
$c\cdot\vol (M)$.
\end{lem}

\begin{proof}
Write $M=\gC\backslash S$ and let $\ga, \gb\in\gC$ be elements which correspond to two
different closed geodesics in $M$ of length $\leq\gep$. Then the axes of $\ga$
and $\gb$ are bounded away from each other, and hence, for large enough $m$, $\langle\ga^m,\gb^m\rangle$ is a non-abelian
free group. It follows that $\{ d_\ga <\gep_s\}\cap\{ d_\gb <\gep_s\}=\emptyset$.

Since $d_\gc (x)\leq \gep+2D_{\{ d_\gc <\gep\}}(x)$, we see that
$$
 \{ D_{\{ d_\ga <\gep\}}<\gep\}\cap\{ D_{\{ d_\gb <\gep\}}<\gep\}=\emptyset.
$$
This implies that if we take, for each connected component $M_{\leq\gep}^0$
of $M_{\leq\gep}$, an $\gep$-ball $B_{\gep}$, whose center lies on the
boundary of $M_{\leq\gep}^0$, then these balls are disjoint and injected.
Thus the number of geodesics of length $\leq\gep$ (which coincides with
the number of connected components of $M_{\leq\gep}$ and hence with the number 
of these ${\gep}$-balls) is $\leq\vol (M)/\vol (B_{\gep})$.
\end{proof}

\begin{proof}[Proof of Theorem \ref{Pre-Hyp}]
Let $\gep =\frac{\gep_s}{3}$, let $M$ be a complete hyperbolic 
$3$-manifold, and let $N\subset M$ be the union of all closed geodesics in $M$
of length $\leq\gep$. Then $N$ consists of 
$\leq\frac{\vol (M)}{\vol (B_{\gep})}$
circles. As in the proof of \ref{thm3}, $\pi_1(M\setminus N)$ is
isomorphic to $\pi_1(\mathcal{R})$ for some 
$\big(\ga,d\cdot\vol(M)\big)$-simplicial complex $\mathcal{R}$.
It follows that $\pi_1(M\setminus N)$ admits a presentation with 
$\leq\frac{\ga d}{2}\vol (M)$ generators, and $\leq \ga d^2\vol (M)$ relations.
Van-Kampen's theorem implies that when adding these circles one by one
to $M\setminus N$, we should add one relation for each (note that each circle
has a neighborhood homeomorphic to a solid torus or a solid Klein bottle) .
Hence, we get a presentation of $\pi_1(M)$
with the same number of generators and with at most 
$\frac{\vol (M)}{\vol (B_\gep)}$ additional relations.
\end{proof}

\begin{rem}\label{55}
In contrast with Theorem \ref{presentation-thm}, Theorem \ref{Pre-Hyp} does 
not yield bounds for the length of the relations. Since for $v$ large enough, 
there are infinitely many complete hyperbolic $3$-manifolds with volume 
$\leq v$, the length of the relations in the above presentations can not be 
bound in terms of $\vol (M)$.
\end{rem}

\begin{rem}
For $S=\BH^2\times\BH^2$ the analog of Theorem \ref{Pre-Hyp} holds. 
The proof is almost the same, except
that in this case, each of the $\leq\frac{\vol (M)}{\vol (B_{\gep /2} )}$ 
connected components of the union $N$ of all closed 
geodesics of length $\leq\gep$ is either a circle or a two dimensional torus
or a Klein bottle.
\end{rem}

However for the symmetric space $S=\SL_3(\BR )/\SO_3(\BR )$ we have only the 
following result which follows directly from Theorem \ref{thm3}.

\begin{prop}
For any $S$, there is a constant $\eta (S)$ such that the fundamental group
of any $S$-manifold $M$ has a generating set of size $\eta (S)\cdot\vol (M)$.
\end{prop}


\section{A quantitative version of Wang's theorem}\label{QuaWang}

Denote by $\rho_S(v)$ the number of isometric classes of irreducible 
$S$-manifolds of volume $\leq v$. 

If Conjecture \ref{conjA} is true then for any
symmetric space, $S$, of non-compact type of dimension $\geq 4$, there is some constant $c=c(S)$,
such that
$$
 \rho_S(v)\leq v^{c\cdot v}
$$
for any $v>0$.
While in dimension $3$, the validity of Conjecture \ref{conjA} would yield an analogous upper bound
for the growth of {\it arithmetic} $3$-manifolds.

Theorems \ref{thmA}(2) implies: 

\begin{thm}\label{QWthm}
Assume that $S$ is neither isometric to 
$\BH^2,\BH^3,\SL_3(\BR )/\SO_3(\BR )$ nor to $\BH^2\times\BH^2$.
Then we have
$$
 \rho_S(v)\leq v^{c(S)\cdot v}.
$$
\end{thm}

And \ref{thmA}(1) implies:

\begin{prop}\label{12.2}
The number of non-compact arithmetic hyperbolic $3$-manifolds of volume
$\leq v$ is at most $v^{cv}$ for some constant $c$.
\end{prop}

\begin{rem}\label{12.3}
Similarly, for $S=\SL_3(\BR )/\SO_3(\BR )$ or $\BH^2\times\BH^2$, the
number of non-compact irreducible manifolds with volume $\leq v$ is bounded
by $v^{cv}$ for some $c$.
\end{rem}

\begin{proof}[Proof of \ref{QWthm}, \ref{12.2} and \ref{12.3}]
Since $\dim (S)\geq 3$, it follows from Mostow's rigidity theorem that an 
irreducible $S$-manifold $M$ is characterized by its fundamental 
group, which, by Theorem \ref{presentation-thm}, has a presentation 
$\pi_1(M)\cong\langle \gS , W\rangle$ with $|\gS |,|W|\leq\eta (S)\vol (M)$ in which all the relations has length $\leq 3$.
A rough estimate of the number
of groups admitting such a presentation yields \ref{QWthm}. 
\end{proof} 

\begin{rem}
It was shown in \cite{BGLM} that for $\BH^n$ when $n\geq 4$ this estimate
is tight. However in the higher rank case it is very likely that $\rho_S(v)$ grows
much slower. It was guessed in \cite{BGLM} that when $\rk (S)\geq 2$
$$
\log\rho_S(v)\approx {c(S)\frac{(\log V)^2}{\log\log V}}.
$$
\end{rem}

\begin{rem}
Since Mostow rigidity does not hold for 
surfaces, our method does not yield a quantitative version for Borel's 
finiteness theorem for arithmetic hyperbolic surfaces of a given genus.
\end{rem}


\section{Some complements}\label{complements}

\subsection{Extending some of the results to non-compact orbifolds} 

In the previous sections we have considered $S$-manifolds of finite volume. It is natural to try 
to generalize the results obtained, to the larger family of $S$-orbifolds of finite volume. 
This amounts to consider general lattices in $G$ instead of just torsion free lattices.

It turns out that some of the main statements could be generalized to  
non-compact $S$-orbifolds, i.e. to general non-uniform lattices 
$\gC\leq G$. We remark that we do not know how to deal with general compact orbifolds.
In the non-compact case, our generalizations rely on the following effective version of Selberg's lemma: 

\begin{lem}\label{i(t.f.)}
There is a constant $i=i(G)\in\BN$ such that any non-uniform arithmetic lattice
$\gC\leq G$ has a torsion free normal subgroup of index $\leq i$.
\end{lem}

\begin{proof}
Let $\mathfrak{g}$ denote the Lie algebra of $G$ and let $n$ be its dimension.
Replacing $G$ by its identity component we can assume it is connected and 
center free, and therefore may be identified with its image under the adjoint 
representation $\textrm{Ad}(G)\leq\textrm{GL}(\mathfrak{g})\cong GL_n(\BR )$. 
As follows from the proof of Margulis' arithmeticity theorem, for any 
non-uniform arithmetic lattice $\gC\leq G$ there is a base $B$ for the
vector space $\mathfrak{g}\cong \BR^n$ with respect to which 
$\gC\leq \textrm{GL}_n(\BQ )$ and is commensurable to $\textrm{GL}_n(\BZ )$.
Replacing this base by a $\BZ$-base for the $\BZ$-span of $\gC\cdot B$ (which
is easily seen to be a $\BZ$-lattice in $\BR^n$), we can assume that $\gC$ is 
contained in $\textrm{GL}_n(\BZ )$.

Let $T\leq\textrm{GL}_n(\BZ )$ be a fixed torsion-free congruence subgroup  
(which exists, for instance, by Selberg's lemma) and let $i=i(G)$ be its index
$$
 i=[\textrm{GL}_n(\BZ ):T].
$$
Clearly, $\gC\cap T$ is torsion free and $[\gC :\gC\cap T]\leq i$.
\end{proof}

The following generalization of Theorem \ref{presentation-thm} follows 
immediately:

\begin{thm}\label{Gpresentation}
There is a constant $\eta (G)$ such that any non-uniform arithmetic 
lattice $\gC\leq G$ has a presentation 
$\gC\cong \langle \Sigma :W\rangle$
with 
$$
 |\Sigma |,|W|\leq\eta (G)\cdot\vol (G/\gC ).
$$
\end{thm}

\begin{proof}
Take a torsion free normal subgroup $\gC_1$ of index $\leq i$ in $\gC$.
Then $\vol (G/\gC_1)=\vol (G/\gC )|\gC/\gC_1|\leq\vol (G/\gC )\cdot i$, and 
$\gC_1$, being torsion free, has a presentation with $\leq\eta'\vol (G/\gC_1)$ generators and 
relations by \ref{presentation-thm}. We should add at most $|\gC/\gC_1|\leq i$ generators and 
$|\gC/\gC_1|^{|\gC/\gC_1|}\leq i^i$ relations to get a presentation for $\gC$.
\end{proof}

However, unlike the case of \ref{presentation-thm}, we do not know how to
bound the lengths of the relations in $W$.

\medskip

The following extends the results of section \ref{QuaWang}.

\begin{thm}\label{countingL}
There is a constant $c=c(G)$ such that for any $v>0$, the number of conjugacy 
classes of non-uniform arithmetic lattices of covolume $\leq v$ is at most
$v^{c\cdot v}$.
\end{thm}

\begin{proof}
We already know (by section \ref{QuaWang}) that for any $v>0$ there are at most 
$v^{c'\cdot v}$ 
conjugacy classes of torsion free lattices of covolume $\leq v$, and therefore 
at most $v^{c''\cdot v}$ conjugacy classes of torsion free lattices of covolume
$\leq i\cdot v$ where $i=i(G)$ is the constant from Lemma \ref{i(t.f.)}.

Let $v_0$ be the minimal covolume of a lattice $\gC\leq G$. Let $\gC\leq G$ be 
a lattice of covolume $\leq v$. By Lemma \ref{i(t.f.)}, $\gC$ contains a 
torsion free normal subgroup $\gC '$ of index $\leq i$. Let $N_G(\gC ')$ be 
the normalizer of $\gC '$ in $G$. Then $N_G(\gC ')$ is a lattice containing 
$\gC$ whose covolume satisfies
$$
 v_0\leq \vol\big( G/N_G(\gC ')\big)\leq v.
$$ 

It follows from Theorem \ref{Gpresentation} that $N_G(\gC ')$ has a generating 
set of size $\leq [\eta\cdot v]$. Thus, for any $j$, the number of subgroups 
of $N_G(\gC ')$ of index $j$ is no more than 
$\big(j+[\eta\cdot v]\big) !^2$ 
(which is a trivial upper bound for the number of index $j$ subgroups of the 
free group of rank $[\eta\cdot v]$). The index $[N_G(\gC '):\gC ]$ is at most 
$v/v_0$, so there are at most
$\big( [v/v_0]+[\eta\cdot v]\big) !^2$ choices for $\gC$ as a
subgroup of index
$[N_G(\gC '):\gC ]$ of $N_G(\gC ')$. Summing over all indexes $\leq v/v_0$, 
we get that there are at most
$[v/v_0]\cdot \big( [v/v_0]+[\eta\cdot v]\big) !^2\leq v^{c'''\cdot v}$ 
choices for $\gC$ as a
subgroup of index $\leq [N_G(\gC '):\gC ]$ of $N_G(\gC ')$. 

Thus, the number of lattices $\gC$ which contain the same $\gC '$ as a 
normal subgroup of index $\leq i$ is at most $v^{c'''\cdot v}$. 
Since there are at most $v^{c''v}$ possibilities for $\gC'$,
it follows that the number
of conjugacy classes of lattices of covolume $\leq v$ is at most 
$v^{c''\cdot v}\cdot v^{c'''\cdot v}\leq v^{c\cdot v}$.
\end{proof}


\subsection{Commensurable growth}
Let us now restrict our attention to a fixed commensurability class. 

\begin{defn}
Two $S$-manifolds $M,N$ are called {\bf commensurable} if they have a common
finite cover. I.e. $\gC_1\backslash S$ is commensurable to $\gC_2\backslash S$ iff $\gC_1$ is
commensurable to some conjugate of $\gC_2$ in $G=\text{Isom}(S)$.
\end{defn}
   
The following definition is natural.

\begin{defn}
The {\bf commensurable growth} $\kappa_M(v)$ of a locally symmetric manifold 
$M$, is the number of non-isometric manifolds commensurable to $M$  
with volume $\leq v$.
The {\bf commensurable growth} $\kappa_\gC (v)$ of a lattice $\gC\leq G$ is the
number of conjugacy classes of lattices commensurable to $\gC$ with
covolume $\leq v$.
\end{defn}


One can define the notion of commensurable growth for arbitrary subgroup
$\gC\leq G$, not necessarily a lattice, as follows: define the 
generalized index between commensurable subgroups
$\gC ,\gC'\leq G$ to be the rational number 
$$
 [\gC :\gC']=\frac{[\gC :\gC\cap\gC']}{[\gC' :\gC\cap\gC']},
$$ 
and use this concept instead of ``covolume'' in the above definition.

Clearly, for a locally symmetric manifold $M=\gC\backslash S$, we have
$\kappa_M(v)\leq\kappa_\gC (v)$. It is natural to ask what is the relation
between these functions. In particular, 
do they have the same asymptotic behavior?    

Another interesting question is what is the relation between the commensurable growth  
and the congruence subgroup problem.

\medskip

We shall now give upper bounds for the commensurable growth of locally 
symmetric manifolds and for its fundamental group when the dimension is $>2$.

\begin{clm}
Let $M=\gC\backslash S$ be an irreducible locally symmetric manifold of dimension $>2$ 
with finite volume. Then there is a constant $c=c(M)$ such that 
$\kappa_M(v)\leq v^{cv}$.
\end{clm}

\begin{proof}
If $M=\gC\backslash S$ is not arithmetic then, by Margulis' criterion for arithmeticity,
the commensurability class of $\gC$ admits a unique maximal element which 
contains all the others, and the result follows by considering the subgroup 
growth of this maximal element. If $M$ is arithmetic, this follows from \ref{fixedO} and from 
\ref{thm1}.
\end{proof}

When $\gC\leq G$ is a non-arithmetic lattice, the above proof applies also to $\kappa_\gC$ and gives 
$\kappa_\gC (v)\leq v^{cv}$ as well. In the arithmetic case, as in 
\ref{countingL}, we obtain similar upper bounds by using the following 
lemma, which can be proved in the same way as Lemma \ref{i(t.f.)}.

\begin{lem}\label{180}
For any commensurability class $\mathfrak{N}$ of arithmetic lattices in $G$,
there is a constant $i=i(\mathfrak{N})$ such that any $\gC\in\mathfrak{N}$ 
contains a torsion free subgroup of index $\leq i$.
\end{lem}

The following is immediate from \ref{fixedO} and \ref{180}.

\begin{lem}
Given a commensurability class $\mathfrak{N}$ of arithmetic lattices in $G$,
there is a constant $\eta =\eta (\mathfrak{N})$ such that any 
$\gC\in\mathfrak{N}$
has a presentation $\gC\cong \langle \Sigma :W\rangle$
with $|\Sigma |,|W|\leq\eta (G)\cdot\vol (G/\gC )$.
\end{lem}

As in the proof of \ref{countingL}, these two lemmas imply:

\begin{prop}
For any lattice $\gC\leq G$, $\kappa_\gC (v)\leq v^{c(\gC )v}$. 
\end{prop}


\subsection{How to construct a simplicial complex for non-arithmetic manifolds}

We shall now explain how to attach simplicial complexes to
non-arithmetic, or more generally to rank-$1$ manifolds, of dimension $\geq 4$. We are not
trying to do it in the most economical way, but just to explain an idea of how this could be done.

The following proposition follows from a rough estimate for the diameter and 
the minimal injectivity radius of compact connected components of the thin 
part.

\begin{prop}\label{rank-1prop}
There are positive constants $\ga =\ga (n)$ and $d=d(n)$ such that 
any compact rank-$1$ locally symmetric manifold $M$ of dimension $n\geq 4$ is 
homotopically equivalent to a 
$\big(d,\ga\cdot\vol (M)^{3n^2+1}\big)$-simplicial complex.
\end{prop}

\begin{proof}
Fix $n$, let $\gep (n)$ be the constant of the Margulis lemma, and let 
$M_{\leq\gep (n)}$ be the thin part of the ordinary $\gep (n)$ thick-thin 
decomposition.

It follows from \cite{BS} (see proposition 3.2 there) that for some constant
$c$, the diameter of any compact connected component of $M_{\leq\gep (n)}$ is
at most $3\log \big( c\cdot\vol (M)\big)$. 
Applying formula 8.5 from \cite{Gr} (page 381), which implies that the 
injectivity radius decreases at most exponentially as one moves along the 
manifold, we get that the injectivity radius at any point belonging to
$M_{\leq\gep (n)}$, and therefore at any point of $M$, 
is at least $c'\cdot\vol (M)^{-3n}$, for some positive constant $c'$. 

The proposition follows by applying a good covering argument with ordinary 
balls of radius $\gep =c'\cdot\vol (M)^{-3n}$. The needed number of balls in such a
cover is $\leq c''\frac{\vol(M)}{\gep^n}\leq\ga\cdot\vol (M)^{3n^2+1}$.
\end{proof}     

In fact, it is easy to obtain stronger results by means of elementary 
computations of volumes of neighborhoods of short closed geodesics. 
We shall demonstrate this in the real hyperbolic case. We remark that in all
other rank-$1$ cases, one can obtain similar estimates by using the same
means, and applying Rauch's comparison theorems (see \cite{CH} 1.10). 
However, also the following estimate is probably not tight, and
it should be possible (and not necessarily very hard) to obtain better 
estimates. 

\begin{prop}
For $n\geq 4$, there are constants $\ga =\ga (n), d=d(n)$, such that any 
compact hyperbolic $n$-manifold $M$ is homotopically 
equivalent to some 
$\big(d,\ga\cdot \vol (M)^{(1+n( [\frac{n-1}{2}]+1))/(n-2-[\frac{n-1}{2}])}
\big)$-simplicial complex. 
\end{prop}

\begin{proof}[Sketched proof for $n=4$]
In this case any connected
component of the thin part of the ordinary thick-thin decomposition (with 
respect to some fixed $\gep$) is a neighborhood of a short closed geodesic 
which is topologically a ball bundle over a circle.
In order to understand the geometry of the thin components it is most 
convenient to look at the upper half-space model for the hyperbolic space
$\BH^4$ (see \cite{BePe} for details). Lift the component so that the short 
closed geodesic is lifted to the line connecting $0$ to $\infty$. 
Then our lifted 
component is a cone, centered by this line. The intersection of this cone with
a horosphere perpendicular to this line is a union of coaxial ellipsoids (with 
respect to the induced $(n-1)$-Euclidean structure on the horosphere). 

Using the fact that any abelian subgroup of $\textrm{SO}_3(\BR )$ 
($\cong$ to the fixator group of this line) is contained 
in a $1$-dimensional torus, one can show, by a simple pigeonhole argument on the powers of the corresponding 
isometry, that if the length of our short closed
geodesic is $a<<\gep$ then, for some constant $c$, the $3$ dimensional
Euclidean ball of radius $c\frac{1}{a^{1/2}}$ is contained in the union of 
the above ellipsoids. This implies that the volume of the component is at least
$c'(1/a^{1/2})^3\cdot a$. Thus $a\geq c''\frac{1}{\vol (M)^2}$. 

Thus the injectivity radius at any point of $M$ is at least 
$\rho =c''\frac{1}{\vol (M)^2}$, 
and one can construct, as above, a simplicial complex with at most
$c'''\frac{\vol (M)}{\rho^4}=\ga\cdot\vol (M)^9$ vertices, all of them of 
degree bounded by some constant $d$.

We remark that the proof for general dimension $n$ uses the fact that any 
abelian subgroup of $\textrm{SO}_{n-1}(\BR )$ is contained in some 
$[\frac{n-1}{2}]$-dimensional torus.
\end{proof}

\begin{rem}
In order to obtain analogous estimates for non-compact rank-$1$ manifolds, one 
should use the thick-thin decomposition with $\gep =\gep \big(\vol (M)\big)$
as above, so that all the components of $M_{\leq\gep}$ would be cusps, and 
then estimate explicitly the function $\gb (\gep)$ which is defined in the
proof of Theorem \ref{thm1},
and detect its influence on the determination of the constant $b=b(\gep )$
in Proposition \ref{simplicial-complex} in order to finally calculate the resulting 
simplicial complex.
\end{rem}


\subsection{Wang's theorem for products of $\SL_2$'s}\label{SL_2}

This paragraph is not precisely a part of the main theme of this paper, but only a part of the 
same subject of mathematics. 
Moreover, we are not presenting any new result here, but only clarify things which are
evidently known to some peoples. The author decided to write this paragraph because it might serve
as a complement to Wang's paper \cite{Wa}.

Wang's theorem states that if $G$ is a connected semisimple
Lie group without compact factors, and $G$ is not locally isomorphic to 
$\textrm{PSL}_2(\BR )$ or $\textrm{PSL}_2(\BC )$,
then for any $v>0$, there are only 
finitely many conjugacy classes of irreducible lattices in $G$ of covolume 
$\leq v$.

In \cite{Wa} Wang didn't consider the case where $G$ is of higher rank and
has factors locally 
isomorphic to $\textrm{PSL}_2(\BR )$ or $\textrm{PSL}_2(\BC )$.
Margulis' arithmeticity theorem implies that if $G$ has both a 
factor which is locally isomorphic to $\textrm{PSL}_2(\BR )$ or 
$\textrm{PSL}_2(\BC )$ and a factor which is not locally isomorphic
$\textrm{PSL}_2(\BR )$ or $\textrm{PSL}_2(\BC )$ then $G$ contains no 
irreducible lattices (see \cite{Mar1} corollary 4.5, page 315). 
It was remarked by Borel (see \cite{Bor} 8.3 and 8.1)
that Wang's argument implies also the finiteness
of the number of conjugacy classes of irreducible lattices for groups locally
isometric to $G_{a,b}=\SL_2(\BR )^a\times\SL_2(\BC )^b$ when 
$(a,b)\neq (1,0),(0,1)$. 
Let as now explain this remark of Borel.

It follows from Margulis' super rigidity theorem that irreducible higher rank
lattices are locally rigid. Thus, the missing 
ingredient in Wang's argument (\cite{Wa}, 8.1) when applied to groups 
locally isometric to $G_{a,b}$ is the following statement (which was also 
noted without proof in \cite{Bor}).

\begin{prop}\label{irr-lim}
Let $G$ be a semi-simple Lie group with no compact factors, and let 
$\Gamma_n\leq G$ be a sequence of irreducible lattices. Assume that 
$(\Gamma_n)$ converges to a lattice $\gD$ in the topology of closed 
subgroups (Hausdorff convergence on compact sets). Then $\gD$ is also 
irreducible. 
\end{prop}

\begin{proof}
Assume that $\gD$ is reducible. Then we can write $G$ as an almost direct 
product $G=G_1\cdot G_2$ in such a way that $\gD$ is commensurable with
$\gD_1\cdot\gD_2$ where $\gD_i=\gD\cap G_i$. Fix a finite generating set for
$\gD$, and for large $n$, denote by $f_n:\gD\to\gC_n$ the homomorphism induced
by sending each generator to the closest element in $\gC_n$. As explained in 
\cite{Wa}, since $\gD$ is finitely presented, $f_n$ is a well defined 
homomorphism whenever $\gC_n$ is close enough to $\gD$.

For any $\gd\in\gD$, $f_n(\gd )\to \gd$.  
We will show that $f_n(\gd )$ is central for each non-central $\gd\in\gD_1$, and for any large enough $n$. Since the center of $G$ is 
discrete, this will imply the desired contradiction. 
     
Fix $\gd\in\gD_1$ non-central. As $\gC_n$ is irreducible, 
we will show that $f_n(\gd )$ is central by showing that its projection
to the second factor $\pi_2(f_n(\gd ))$ is the unit element in $G_2$. 

$\gD_2$ is a lattice in $G_2$.
Let $\{\gd_{2,1},\gd_{2,2},...,\gd_{2,k}\}\in\gD_2$ be a finite set of 
generators for $\gD_2$. By Borel's density theorem, 
$\{\text{Ad}(\gd_{2,1}),...,\text{Ad}(\gd_{2,k})\}$ generates the algebra 
$$
 \langle \text{Ad} (G_2)\rangle \leq
 \textrm{End}(\mathfrak{g}_2)
$$ 
(here $\mathfrak{g}_2$ denotes the Lie algebra 
of $G_2$). Since this algebra is finite dimensional, it is generated by
$\{ \text{Ad}\big(\pi_2 ( f_n (\gd_{2,i} ))\big)\}_{i=1}^k$ 
whenever $n$ is 
large enough. Since $G_2$ is semi-simple, the adjoint representation 
$\text{Ad}:G_2\to\textrm{GL}(\mathfrak{g}_2)$ has no invariant vectors.

Let $\gep_n=\pi_2 (f_n (\gd ))$. Since $f_n (\gd )$ is close to $\gd$,
$\gep_n=\pi_2 (f_n (\gd ))$ is close to the identity of $G_2$. We can therefore
assume that $\gep_n$ is contained in an identity neighborhood of $G_2$ where 
$\log =\exp^{-1}:G_2\to\mathfrak{g}_2$ is a well defined diffeomorphism.
As $\gd$ commutes with each $\gd_{2,i}$, $\gep_n$ commutes with each 
$\pi_2 (f_n (\gd_{2,i} ))$, and it follows that
$$
 \text{Ad}\big(\pi_2 ( f_n (\gd_{2,i} ))\big)(\log \gep_n )=
 \log \gep_n,
$$
which in turn implies $\log \gep_n=0$, i.e. $\gep_n=1$.
\end{proof}

Together with Proposition \ref{irr-lim}, the original argument from \cite{Wa} 8.1 gives:

\begin{thm}[Wang's theorem]\label{WT}
Let G be a connected semi-simple Lie group without compact factors, which is 
not locally isomorphic to $\SL_2(\BR )$ or $\SL_2(\BC )$. Then for any $v>0$
there are only finitely many conjugacy classes of irreducible lattices in $G$
with covolume $\leq v$. 
\end{thm}  
       
\begin{rem}
In \cite{BP}, Borel and Prasad established a very strong and general finiteness result, but they
omitted the cases of $G=G_{a,b}$ by requiring absolute rank $\geq 2$. This requirement was used in 
their proof of the stronger finiteness statement, in which the ambient group $G$ can be varied. 
However, Prasad remarked to the author that, when $G$ is fixed, this requirement is unnecessary
in their argument, and hence, the finiteness of the number of 
conjugacy classes of arithmetic lattices of covolume $\leq v$ in $G$
could be proved also by using their methods.
\end{rem}

\begin{rem}
More generally, for any $G$, the finiteness statement holds for lattices (not necessarily 
arithmetic) which are irreducible with respect to the $\SL_2$ factors of $G$.
I.e. for the set of conjugacy classes of lattices in $G$ which project densely
to any factor of $G$ which is locally isomorphic to $\SL_2(\BR )$ or
$\SL_2(\BC )$.
\end{rem}

\begin{Ack}
I would like to thank 
Shahar Mozes (my Ph.D. adviser) for his guidance throughout 
the last few years and for insightful suggestions for this work,
to Pierre Pansu for hours of helpful conversations and for 
suggestions which were essential to this research,
to Alex Lubotzky for many discussions, suggestions and ideas, 
and to Uri Bader, Yair Glasner, and Yehuda 
Shalom for many discussions and clever suggestions. Proofs of some of the 
statements presented here were established during these discussions. 
I would also like to thank Emmanuel Breuillard, Assaf Naor and Gopal Prasad 
for some remarks concerning early versions of this paper.
Finally, I would like to express my sincere gratitude to the anonymous referees 
for their careful reading of the manuscript and for many 
suggestions and corrections which tremendously improved the exposition of this paper.

\end{Ack}

\end{document}